\documentclass[12pt,a4]{amsart}
\usepackage{amscd,amsmath,amssymb,amsfonts,ifthen,hyperref,tikz-cd} 
\usepackage{calligra}
\usepackage{mathtools}

\usepackage{bbold,tikz}
\usepackage{color} 
\usepackage{graphicx}
\usepackage[all]{xypic}

 \usepackage{times}
\usepackage{fourier}

\usepackage[T1]{fontenc}
\newtheorem{theorem}{Theorem}
\newtheorem{proposition}[theorem]{Proposition}
\newtheorem{claim}{Claim}[subsection]
\newtheorem{thm}[claim]{Theorem}
\newtheorem{lem}[claim]{Lemma}
\newtheorem{cor}[claim]{Corollary}
\newtheorem{example}[claim]{Example}
\newtheorem{prop}[claim]{Proposition}
\theoremstyle{definition}
\newtheorem{defn}[claim]{Definition}

\newtheorem{rem}[claim]{Remark}

\usepackage{geometry}
 \author[VQ Bao]{V\~o Qu\^oc Bao}
 \email[VQ Bao]{vqbao@math.ac.vn}
 \author[PH Hai]{Ph\`ung H\^o Hai}
 \email[PH Hai]{phung@math.ac.vn}
 \author[DV Thinh]{D\`ao Van Thinh}
 \email[DV Thinh]{dvthinh@math.ac.vn}
 \address[VQ Bao, PH Hai, DV Thinh]{Institute of Mathematics, Vietnam Academy of Science and Technology}  

\title[Cohomology of the stratified fundamental group]{Cohomology of the stratified fundamental group of curves}

\makeatletter
\@namedef{subjclassname@2020}{\textup{2020} Mathematics Subject Classification}
\makeatother

\begin{document}
\maketitle
\begin{abstract}
	Let \(X\) be a smooth projective curve over an algebraically closed field of positive characteristic. Using the Frobenius-divided description of stratified bundles, we compare the cohomology of  stratified bundles on \(X\) with the cohomology of the stratified fundamental group. We prove that the natural comparison morphisms are isomorphisms in all degrees, yielding a positive-characteristic analogue of the de Rham  \(K(\pi,1)\)-scheme. In particular, the higher cohomology groups vanish in degrees at least two. As a corollary, we recover the classical result, due to Serre and Pacheco-Stevenson, that the étale fundamental group has p-cohomological dimension one. We also investigate the behavior of these comparison morphisms for ind-stratified bundles. For elliptic curves the comparison holds for all ind-stratified bundles in the supersingular case, and for isotypical ind-stratified bundles in the ordinary case; we give an ordinary elliptic counterexample for an infinite sum of rank-one objects. In higher genus, by constructing explicit Frobenius-invariant ind-bundles and introducing the notion of Cartier nilpotency order, we exhibit a different failure mechanism already inside the unipotent isotypical component.
\end{abstract}


\section{Introduction}
Let $X$ be a smooth geometrically connected algebraic variety over an algebraically
closed field $k$ of characteristic $p>0$. 
A \textit{stratified sheaf} is a quasi-coherent $\mathcal O_X$-module equipped with an action of the sheaf
of differential operators $\mathcal D (X/k)$. The terminology ``stratified'' is due to Grothendieck,
in connection with the infinitesimal cohomology,
see, e.g., \cite[Section~2]{BO78} for more details. It is shown that a stratified sheaf which is
coherent as an $\mathcal O_X$-module is automatically locally free, hence will be referred to as a
``stratified bundle''. As a consequence, the category of stratified bundles is an abelian, rigid tensor
category.

The category of stratified sheaves is denoted by
$\mathrm{str}(X/k)$, the full subcategory  of stratified bundles is denoted by
$\mathrm{str}^\mathrm{coh}(X/k)$. The cohomology of stratified sheaves as the right derived
functors of the hom-set functor in $\mathrm{str}(X/k)$ was studied by Ogus in \cite{Og75}, see
also \cite{Hai13}.

An action of $\mathcal D (X/k)$ on a sheaf $\mathcal E$ 
determines a descending chain of subsheaves $\mathcal E^{(n)}$ which are obtained
from one another by the Frobenius pull-back: $\mathcal E^{(n)}\cong
F^*\mathcal E^{(n+1)}$. Such a sequence $(\mathcal E^{(n)})$ together with
the isomorphisms $\sigma_n:\mathcal E^{(n)}\longrightarrow
F^*\mathcal E^{(n+1)}$ gives an alternative definition
for a stratified sheaf, in case $\mathcal E$ is coherent, 
\((\mathcal E^{(n)},\sigma_n)\) is called a \textit{flat bundle}, cf.  \cite[Section~1]{Gie75}. 
Some other authors also call it an \textit{F-divided  bundle} to generalize it to a more general
setting, see, e.g., \cite{BHdS23}.

Let $x$ be a $k$-point of $X$. The fiber functor at $x:$  
\begin{align*}
 x^*:   \mathrm{str}^\mathrm{coh}(X/k) \longrightarrow \mathrm{Vec}_k,\quad 
    (\mathcal{V}, \nabla) \mapsto \mathcal{V}_{x},
\end{align*}
is faithfully exact (since the kernel and cokernel of a map in $\mathrm{str}^\mathrm{coh}(X)$
are always locally free). Hence,
Tannakian duality (cf. \cite[Theorem 2.11]{DM82})  yields an affine group
scheme, denoted by $\pi(X,x)$, with the property that $x^*$ gives an equivalence:
\begin{align}\label{eq_tannakian_equiv}
 x^*:   \mathrm{str}^\mathrm{coh}(X/k) \overset{\cong}{\longrightarrow} 
 \mathrm{Rep}^\mathrm{f}(\pi(X,x)).
\end{align} 
This equivalence extends to an equivalence between the ind-category 
$\mathrm{str}^\mathrm{ind}(X/k)$ and the category $\mathrm{Rep}(\pi(X,x))$.
Here  $\mathrm{str}^\mathrm{ind}(X/k)$ denotes the category of stratified sheaves
which can be presented as inductive limits of stratified bundles. We shall refer to
them as \textit{ind-stratified bundles}. The group $\pi(X,x)$ can be seen
as an analog in positive characteristic of the differential fundamental group. It was 
implicitly studied in the work  \cite{Gie75} of Gieseker and more recently in the work
of dos Santos \cite{dSa07} and subsequently by many authors, see, e.g., 
\cite{EM10, Kin15,ESB16,vdP19, Sun19}. 

In what follows we shall fix the point $x$
and abbreviate $\pi(X,x)$ to $\pi(X)$ for simplicity. 

A central theme of this work is the interaction between the  cohomology theory of stratified sheaves and the rational cohomology attached to the stratified fundamental group. Although Tannakian duality provides an equivalence between stratified bundles and finite dimensional representations of the stratified fundamental group, the corresponding cohomology theories are a priori defined in rather different ambient categories. In particular, the stratified cohomology is computed inside the larger category of all stratified sheaves and is governed by inverse systems arising from Frobenius descent, while the group cohomology is defined through injective resolutions in the representation category of an affine group scheme. Understanding the precise relationship between these two theories therefore requires a careful analysis of how Frobenius-divided structures behave under limits and extensions. More precisely, the equivalence in \eqref{eq_tannakian_equiv} 
implies that the group cohomology of $\pi(X)$ with coefficients in a representation $V$ is the same
as the cohomology of the corresponding stratified bundle $\mathcal V$ computed in the category 
$\mathrm{str}^\mathrm{ind}(X/k)$ (that is the derived functors of the hom-set functor computed in
$\mathrm{str}^\mathrm{ind}(X/k))$. 
Consequently, there is a natural map
\begin{equation}\label{eq_comparison}
\delta^i:\mathrm{H}^i(\pi(X),V)\longrightarrow 
\mathrm{H}^i_\mathrm{str}(X,\mathcal V),
\end{equation}
where $\mathcal V$ is a stratified bundle and $V$ its fiber at $x$, since the functors on both sides of
this natural map are derived functors of the hom-set functors in the corresponding categories. In
particular $\delta^0$ is an isomorphism. Further, as extensions of vector bundles are again vector bundles,
$\delta^1$ is also an isomorphism.

However, $\mathrm{str}^\mathrm{ind}(X/k)$
is much smaller than
the ambient category $\mathrm{str}(X/k)$. Therefore, the maps $\delta^i$ are not expected
to be injective or surjective. Our first main finding is that, for curves, 
$\delta^i$, $i\geq 0$, are nevertheless bijective. 
\begin{theorem}[Theorem \ref{thm_delta_stratified}]
For a smooth projective curve $X$ over $k$ the comparison maps $\delta^i$ in \eqref{eq_comparison} are isomorphisms for all $i\geq 0$
and all stratified bundles $\mathcal V$. As a consequence, the cohomology of $\pi(X)$ vanishes from
degree $i\geq 2$.
\end{theorem}

For smooth, complete varieties over an algebraically closed field, 
the \'etale fundamental group is the profinite quotient of the stratified fundamental 
group.  As a consequence, we recover the following results :
\begin{enumerate}
\item (Proposition \ref{prop_compare_etale}) 
The restriction map $\rho^i:\mathrm{H}^i(\pi^\text{\'et}(X),V)\longrightarrow 
\mathrm{H}^i(\pi(X),V)$ given by the canonical surjective map 
$\pi(X)\longrightarrow \pi^\text{\'et}(X)$, is
bijective for any $i\geq 0$ and any representation $V$ of $\pi^\text{\'et}(X)$.
As a consequence, the cohomology of $\pi^\text{\'et}(X)$ vanishes from degree $i\geq 2$.
\item (Corollary
\ref{cor_nilpotent_cohomology}, Proposition \ref{thm:pi-uni-free})
Let $\pi^\mathrm{uni}(X)$ be the Tannakian dual to the full subcategory of nilpotent
stratified bundles. Then $\mathrm H^i(\pi^\mathrm{uni}(X),V)\cong\mathrm  H^i(\pi(X),V)$
As a consequence, the cohomology of $\pi^\mathrm{uni}(X)$ vanishes from degree $i\geq 2$
so that $\pi^\mathrm{uni}(X)$ is a free pro-\'etale, pro-unipotent group scheme.
\end{enumerate}
The property (1) is due to Serre \cite{Ser90} and Pacheco-Stevenson \cite{PS00}.

 From a conceptual point of view, our comparison theorem suggests that smooth projective curves in positive characteristic should be regarded as ``stratified \(K(\pi,1)\)-schemes''
(cf. \cite{BHT25,HBB26, BN26} and references therein). Namely, the entire stratified cohomology of coherent stratified bundles is recovered from the group cohomology of the stratified fundamental group. This phenomenon is reminiscent of the classical topological situation where the cohomology of local systems on a \(K(\pi,1)\)-space is completely determined by the cohomology of its fundamental group.

We note an important difference: the stratified cohomology groups vanish already in degrees at least \(2\), that is, stratified cohomology does
not satisfy the property of a Weil cohomology.
This differs from the situation for the differential fundamental group in characteristic \(0\) (cf.
\cite{BHT25}) where the cohomology completely vanishes only from degree larger than 2. 
Such a phenomenon was already known in the work of Gieseker, who showed the ``rigidity'' 
of regular singular stratified bundles in the formal neighborhood of a point 
\cite[Theorem~3.3]{Gie75}.

Motivated by the non-vanishing of the first de Rham cohomology of smooth projective
curves of genus $g>0$ \cite{BHT25}, we ask a similar question for stratified cohomology
of curves in positive characteristic. It turns out that, unlike the de Rham cohomology for curves
in characteristic 0, in positive characteristic, the first stratified cohomology group may vanish.
To see that, we carry out some computations for finite stratified bundles, that is \'etale trivializable vector bundles. 
Our explicit computation on the curves in Section \ref{sect_Cartier_nilpotency}
also shows that the non-vanishing of stratified cohomology in degree 1 depends on some
invariants of the curves. This point of view suggests many interesting questions, which
we hope to address in the future. Our second main finding is the following.
\begin{proposition}[Corollary \ref{cor_non-vanishing} and Appendix \ref{proof_non-example}]
\begin{enumerate}
\item 
Let $X$ be an ordinary curve of genus larger than $1$. 
Let $\mathcal E$ be an ordinarily \'etale trivializable vector bundle. 
Then $\mathrm h^1_\mathrm{str}(X, \mathcal E) \neq 0.$
\item There exists an \'etale trivializable vector bundle on a hyperelliptic curve for which the 
first stratified cohomology is equal to zero. 
\end{enumerate}
\end{proposition}

Another important aspect of the paper is that the comparison theorem may fail for ind-stratified bundles. Already on an ordinary elliptic curve, an infinite sum of rank-one stratified bundles can produce a non-zero derived inverse limit and hence a failure of comparison; see Example \ref{example_ordinary_elliptic_ind}. For curves of higher genus we exhibit a different phenomenon. We construct explicit Frobenius invariant iterated extensions on hyperelliptic curves (see Section \ref{sect_comparison_ind}) and introduce the notion of Cartier nilpotency order (cf. Proposition \ref{Prop_Cartier_nilpotency}). The nilpotency order can grow unboundedly along a sequence of extensions, leading to an ind-stratified bundle for which the natural comparison map is no longer bijective even though the whole construction stays in the unipotent isotypical component. These examples highlight the essential difference between coherent and ind-coherent phenomena. They also emphasize that, while group cohomology commutes with inductive limits, stratified cohomology is computed in terms of inverse limit functors (see Section \ref{sect_str_cohomology}), which need not commute with direct limits. 

Our main finding in this direction is the following.
\begin{proposition}[Proposition \ref{prop-non-example} and Lemma \ref{lem_C_E3}]
Let $X$ be a hyperelliptic curve of genus 2 and Hasse-Witt rank 1 over an algebraically closed field $k$ of characteristic $3$. There exists an ind-stratified
bundle $\mathcal E$, which is an infinite extension of the trivial stratified bundle, such that
the canonical map
$$\delta^1: \mathrm{H}^1(\pi(X),E)\longrightarrow \mathrm{H}^1_{\mathrm{str}}(X,\mathcal E)$$
is not bijective. In fact, the target space is infinite dimensional while the source space is zero.\end{proposition}

Let us briefly describe the structure of our work. In Section \ref{sect_str_bundles},
we recall the definitions of stratified bundles and the stratified cohomology following
Gieseker \cite{Gie75} and Ogus \cite{Og75}. 
Each stratified bundle $\mathcal E$ can be equivalently given as a sequence of 
vector bundles $\mathcal E^{(n)}$ such that $F^*\mathcal E^{(n+1)}\cong \mathcal E^{(n)}$,
$F$ being the absolute Frobenius map.
The $i$-th stratified cohomology is computed from the inverse limit of the pro-system
$\mathrm{H}^i(X,\mathcal E^{(n)})$, in which the transfer maps are given by the
action of the Frobenius. 
Thus, the computation of stratified cohomology is reduced to the
study of the Frobenius action on the cohomology groups of sequence of vector bundles
related to each other by the Frobenius pull-back.
We show the finiteness and vanishing properties of this cohomology, cf. 
Proposition \ref{prop_finiteness}.  

In Section \ref{sect_comparison}, we define a natural map from the group cohomology of
the stratified fundamental group to the stratified cohomology.
We show that for smooth projective curves the maps are bijective, cf. Theorem \ref{thm_delta_stratified}. 
Since stratified cohomology on projective curves vanishes
in degree larger than 1, we obtain a similar result for the fundamental group cohomology.
This vanishing theorem also allows us to compare the cohomology of the stratified 
fundamental group with the cohomology of the \'etale fundamental group and its
pro-unipotent quotient in $k$-linear representations (Proposition \ref{prop_compare_etale} and 
Corollary \ref{cor_nilpotent_cohomology}). A key ingredient here is the fact, due to 
dos Santos \cite{dSa07}, that any unipotent quotient of the stratified 
fundamental group is finite. 
In Subsection \ref{sect_etale_trivial} we compare the stratified  
cohomology of finite stratified bundles with the \'etale cohomology of their
$p$-torsion sub-sheaves of Frobenius invariant sections. Then we compare
the stratified cohomology of finite stratified bundles with cohomology of the corresponding
$p$-torsion \'etale sheaves. Our motivation is to study the non-vanishing of the first
cohomology group. 

In Section \ref{sect_comparison_ind}, we consider the case of ind-stratified bundles. For a supersingular elliptic curve the comparison maps $\delta^i$ remain bijective for arbitrary ind-stratified bundles. For an ordinary elliptic curve the same holds for isotypical ind-stratified bundles and their finite direct sums, but it can fail for infinite sums of isotypical components; we give a short explicit counterexample. Section \ref{sect_Cartier_nilpotency} introduces the notion of Cartier nilpotency of Frobenius invariant vector bundles. We construct, on a hyperelliptic curve of genus $g=2$, a sequence $\mathcal E_n$ of iterated extensions of the trivial bundles which are Frobenius invariant. Under the abstract Cartier condition \eqref{condition_Cartier} at $\mathcal E_3$, we show that the order of Cartier nilpotency of $\mathcal E_n$ tends to infinity with $n$. Section \ref{sect_genus_2} verifies this condition on explicit curves in characteristic three. As a consequence, the first stratified cohomology of $\mathcal E=\varinjlim \mathcal E_n$ is infinite dimensional and hence the comparison map $\delta^1$ is not bijective. In contrast with the ordinary elliptic example, this second failure occurs entirely inside the unipotent isotypical component. 

%
%
%

\section{Stratified bundles and stratified cohomology}
\label{sect_str_bundles}

\subsection{Stratified bundles} 
Let $k$ be an algebraically closed field of characteristic $p>0$ 
and let $X/k$ be a smooth, geometrically connected  scheme. 
\subsubsection{} \label{sect_2.1.1}
Let $\mathcal{D}(X/k)$ be  the sheaf of differential 
operators on $X/k$, cf. \cite[IV.16.8]{EGA4}.
A stratified sheaf on $X/k$ is a quasi-coherent $\mathcal O_X$-module $\mathcal V$ equipped 
with an action of $\mathcal D (X/k)$, that is, a $k$-linear map of sheaves of algebras
$$\nabla:\mathcal D (X/k)\longrightarrow \mathcal{E}\textit{nd}_k(\mathcal V).$$
$\nabla$ is also called a stratification on $\mathcal V$. 
The category of stratified sheaves over $X/k$ is denoted by $\mathrm{str}(X/k).$
We denote by $\mathrm{str}^\mathrm{coh}(X/k)$ the full subcategory of $\mathcal O_X$-coherent
stratified sheaves. According to a theorem of Grothendieck, an $\mathcal O_X$-coherent
module equipped with a stratification is always locally free \cite[Prop.~2.16]{BO78} or \cite[Lemma~6]{dSa07}.  
We shall therefore address objects of $\mathrm{str}^\mathrm{coh}(X/k)$ as stratified bundles. 

Thus $\mathrm{str}^\mathrm{coh}(X/k)$ is a rigid tensor category
with the unit object being the trivial connections $(\mathcal O_X,d)$
where $d$ denotes the natural inclusion of $\mathcal D (X/k)$ into
$\mathcal E\textit{nd}_k(\mathcal O_X).$
\subsubsection{} \label{sect_stratified_hom}
Let $(\mathcal E,\nabla)$ be a stratified sheaf on $X/k.$ We define
$\mathcal E^{(i)}$ the subsheaf of sections of $\mathcal E$ 
annihilated by differential 
operators of order less than $p^i$ and vanishing on constant functions.
In other words, let $\mathcal{D}^{<p^{i}}$ be the subsheaf of differential
operators of order less than $p^{i}$. Then  
\begin{align*}
\mathcal E^{(i)} &  = \mathcal{H}\textit{om}_{\mathcal{D}^{<p^{i}}}(\mathcal{O}_X, \mathcal E),
\end{align*}
in particular, $\mathcal E^{(0)}=\mathcal E$. 
We define the sheaf of horizontal sections of 
$(\mathcal E,\nabla)$ to be
$$\mathcal E^{\nabla}:= \bigcap_i \mathcal E^{(i)} = \varprojlim_i \mathcal E^{(i)}.$$
We have
$$\mathrm{Hom}_\text{str}((\mathcal O_X,d),(\mathcal E,\nabla))=
\mathrm{H}^0(X,\mathcal E^{\nabla})=\bigcap_i\mathrm{H}^0(X, \mathcal E^{(i)})= 
\varprojlim_i \mathrm{H}^0(X,\mathcal  E^{(i)}),$$
here $\mathrm{Hom}_\text{str}$ is the abbreviation of $\mathrm{Hom}_\text{str$(X/k)$}$.
This space is always finite dimensional, cf. \cite[Lemma~1.5]{Hai13} or \ref{sect_2.1.4} below. 

\subsubsection{} \label{sect_2.1.3}
According to a result of Cartier, the Frobenius pull-back of $\mathcal E^{(1)}$ is
isomorphic to $\mathcal E$, that is there is an isomorphism
$$\sigma_0:F^*\mathcal E^{(1)}\cong \mathcal E,$$
where $F: X\longrightarrow X$ is the absolute Frobenius map,
which on the structure sheaf is given by 
$F:\mathcal O_{X}\longrightarrow \mathcal O_X$,
$s\mapsto s^p.$
There is a $p$-linear map from $\mathcal E^{(1)}$ to $F^*\mathcal E^{(1)}$.
Composing it with $\sigma_0$ we obtain a $p$-linear map $\psi_1:\mathcal E^{(1)}
\longrightarrow \mathcal E$.

More generally, the Frobenius pull-back of 
$\mathcal E^{(i+1)}$ is isomorphic to $\mathcal E^{(i)}.$ 
As a consequence, there is an alternative
definition of  stratified sheaves as infinite sequences of sheaves
$(\mathcal E^{(i)})_{i\geq 0}$ together with isomorphisms
$\sigma_i: F^*\mathcal E^{(i+1)}\cong \mathcal E^{(i)}$, for all $i\geq 0$, cf. \cite{Gie75}.
This data is called a flat bundle by Gieseker.
We shall denote by $\psi_n$ the resulting $p^n$-linear map $\mathcal E^{(n)}\longrightarrow
\mathcal E$. 
 We notice that, on a projective smooth variety, the stratified bundle
associated to a flat bundle is (up to an isomorphism) independent of the choice
of the isomorphism $\sigma_i$ \cite[Lemma~1.7]{Gie75}.


\subsubsection{} 
We notice that $\mathrm{str}^\text{ind}(X/k)$ is strictly smaller than $\mathrm{str}(X/k)$. For instance, consider $\mathcal D(X/k)$ as a stratified sheaf on itself, then it cannot be presented as the union of its coherent
stratified subsheaves, since any proper subsheaf cannot contain the unit section. 
\subsection{Stratified cohomology} \label{sect_str_cohomology}
When it is clear which stratification is referred to, we shall omit the symbol for it in the
notation of a stratified sheaf, thus we shall write simply $\mathcal E$ instead of
$(\mathcal E,\nabla)$. 
The $0$-th stratified cohomology group of a stratified sheaf $\mathcal E$ 
is defined to be
$$\mathrm H^0_\text{str}(X/k, \mathcal E) :=
\mathrm{Hom}_\text{str}(\mathcal O,\mathcal E)=
 \mathrm{H}^0(X,\mathcal{E}^{\nabla}),$$
 where $\mathcal O$ is equipped with the trivial stratification \ref{sect_2.1.1}. 
This is a left exact functor from $\text{str}(X/k)$ to the category of
$k$-vector spaces. We notice that $\text{str}(X/k)$ has enough injectives.
Hence, we can define the
derived functors $\mathrm H^i_\text{str}(X/k,-)$, which are the higher
stratified cohomology functors. In particular, we have
$$\mathrm H^i_\text{str}(X/k,-)\cong \mathrm{Ext}^i_\mathrm{str}(\mathcal O_X,-).$$

\subsubsection{Computation} \label{sect_Ogus}
Following Ogus (cf. \cite{Og75}), 
we provide some computations of the stratified cohomology
in terms of the inverse limits. 
As in \ref{sect_stratified_hom}, the functor $\mathrm H^0_\text{str}(X/k,-)$ can be seen as
the composed functor
$$\mathrm H^0_\text{str}(X/k,\mathcal E)=\varprojlim_n \mathrm H^0(X, \mathcal E^{(n)}).$$
Therefore, there is a spectral sequence abutting to it. As a consequence, we
have exact sequences (cf. \cite[Theorem 2.4]{Og75})
for any $i \geq 1$:
$$0 \rightarrow {\varprojlim_n}^1\text H^{i-1}(X, \mathcal E^{(n)}) 
\rightarrow\text H^i_\text{str}(X/k, (\mathcal E,\nabla)) \rightarrow 
\varprojlim_n \text H^i(X, \mathcal E^{(n)}) \rightarrow 0.$$
\subsubsection{} \label{varprojlim1}
Recall that the limit on inverse systems of vector spaces 
(or modules) is a left exact functor. Its derived functors, except
possibly for the first one, vanish. Further, for an inverse system 
$(X_i,t_i: X_{i+1}\longrightarrow X_i)_{i\geq 0}$
of vector spaces (or modules over a ring), we have an exact sequence
$$0\longrightarrow \varprojlim X_i\longrightarrow\prod_iX_i
\stackrel d\longrightarrow \prod_iX_i\longrightarrow 
{\varprojlim}^1X_i\longrightarrow 0,$$
where the map $d$ sends a tuple $(x_i)$ to the tuple $(x_i-t_i(x_{i+1}))$.
As a consequence, $\varprojlim{}^1(X_i) = 0$ if and only if the map $d$ is surjective. And we have the following sufficient condition for that.
\begin{prop} 
If $X_i$ is complete and separated with respect to the filtration
\[
F^k(X_i)=\operatorname{Im}(X_{i+k}\longrightarrow X_i)
\]
for every $i$, then $d$ is surjective.
\end{prop}
\begin{proof}
For $j\geq i$, write
\[
t_{i,j}:=t_i\circ t_{i+1}\circ\cdots\circ t_j:X_{j+1}\longrightarrow X_i.
\]
Given $y=(y_i)_i\in\prod_iX_i$, consider the partial sums
\[
x_i^{(N)}:=y_i+\sum_{j=i}^{N}t_{i,j}(y_{j+1})\qquad(N\geq i).
\]
Their tails lie in arbitrarily deep terms of the filtration on $X_i$.
Completeness and separatedness therefore give a well-defined limit
$x_i=\lim_Nx_i^{(N)}\in X_i$. Passing to the limit in
$x_i^{(N)}-t_i(x_{i+1}^{(N)})=y_i$ (for $N\geq i+1$) gives $d(x)=y$.

The preimage $x$ need not be unique: any two preimages differ by an element
of $\ker(d)=\varprojlim_iX_i$. Separatedness is used for uniqueness of the
limit of the partial sums, rather than for uniqueness of a preimage under $d$.
\end{proof}
In particular, if the spaces $X_i$ are finite dimensional, then for each $i$
the images of $X_j\to X_i$ stabilize as $j\to\infty$. Thus the system satisfies
the Mittag--Leffler condition and its first derived limit vanishes. Moreover,
if the dimensions of the spaces $X_i$ are uniformly bounded (i.e.,
$\mathrm{dim}(X_i)$ is bounded above by a constant $c$ for all $i$), then
$\varprojlim X_i$ is finite-dimensional.

\begin{prop}\label{prop_finiteness}
Let $X$ be a smooth connected projective variety over $k,$ and 
$\mathcal E = \{\mathcal E^{(n)} \}$ be a stratified bundle. Then 
the stratified cohomology group can be computed as
the inverse limit of the Zariski cohomology of $\mathcal E^{(n)}$:
$$\mathrm H^i_\mathrm{str}(X/k,\mathcal  E)=
\varprojlim_n \mathrm H^i(X, \mathcal E^{(n)}).$$
Thus $\mathrm H^i_\mathrm{str}(X/k,\mathcal  E)$ are finite dimensional and vanish for
$i>\dim X$.  
\end{prop}
\begin{proof}
We show that the cohomology groups $\mathrm H^i_\mathrm{str}(X/k,\mathcal  E)$ are
finite dimensional for all $i$ and $\mathcal E$. For this, we first recall some 
properties of the Hilbert polynomial of $\mathcal E.$ 
Let $\mathcal{O}_X(1)$ be an ample line bundle. For a bundle $B$ of rank $r > 0,$ 
set $B(m) = B \otimes_{\mathcal{O}_X} \mathcal{O}_X(m)$, and
define the Hilbert polynomial of $B$ relative to $\mathcal{O}_X(1)$ to be
$$p_B(m) = \frac{1}{r} \chi(X,B(m)) \in \mathbb{Q}[m] \quad \text{ for $m\gg 0,$}$$
Additionally, we define the slope of $B$ to be $\mu(B) = \deg(B)/r.$ 
$B$ is said to be semistable (or $\mu$-semistable) if for all subsheaves $U \subset B$ 
we have $\mu(U) \leq \mu(B).$ 

Consider the stratified bundle $\mathcal E = (\mathcal E^{(n)},\sigma_n)_{n \in \mathbb{N}}, \sigma_n: F^*\mathcal E^{(n+1)} \cong \mathcal E^{(n)},$ of rank $r >0.$ Since 
$$\mu((F^n)^*\mathcal{G}) = p^n \mu(\mathcal G),$$ we see that $\mathcal E^{(n)}$ is $\mu$-semistable when $n$ is large enough. Moreover, by \cite[Corollary 2.2]{EM10}, the Hilbert polynomials of $\mathcal E^{(n)}$ are the same and equal to that of $\mathcal O_X.$ Now Langer's boundedness theorem \cite[Theorem 4.2]{La04} tells us that the family $\{ \mathcal E^{(n)}\}$ is bounded, i.e., there exists a scheme $S$ of finite type over $k$ and an $S$-flat family $\mathcal F$ of coherent sheaves on the fibers of $X \times S \rightarrow S,$ which contains $\{ \mathcal E^{(n)} \}.$ 
From the above argument, we have that the family of bundles 
$\{\mathcal E^{(n)} \}_{n \in \mathbb{N}}$ comes from a $S$-flat coherent sheaf, 
say $\mathcal F,$ over $X \times S,$ where $S$ is a scheme of finite type over $k$. 
Hence, to show the uniformly boundedness of $\text H^i(X, \mathcal E^{(n)}),$ 
it is enough to consider the case $S$ is affine. In that case, set $S = \mathrm{Spec }(A),$ 
for some finitely generated $k$-algebra $A.$ Since $X$ is proper, the projection map 
$X \times S \rightarrow S$ is proper too. Hence, by (the Grothendieck complex), 
there is a finite complex $K^\bullet$ of finitely generated $A$-modules such that 
$$\mathrm H^i(K^\bullet \otimes_A B) = 
\mathrm H^i((X\times \mathrm{Spec }(B), \mathcal F \otimes_A B),$$
for any $A$-algebra $B.$ So we deduce our first claim by taking $B$ to be $k(s)$ for any point $s \in S.$ 
\end{proof}

\begin{rem} 
\begin{itemize}
\item[i)] We note the fundamental difference between stratified cohomology in positive
characteristic and its characteristic-zero counterpart, namely de Rham
cohomology. While the de Rham cohomology vanishes only in degrees larger than two times the 
dimension of the variety, stratified cohomology vanishes already in degrees larger than
the dimension of the variety.
\item[ii)] In a recent preprint \cite{Xia25}, the author also proves the finiteness of the stratified cohomology on a smooth and proper scheme by using a similar method. Since he does not assume the projectivity of $X,$ the proof there is more technical. 
\end{itemize}
 \end{rem}

\section{Comparison with fundamental group cohomology}
\label{sect_comparison}
\subsection{The stratified fundamental group and cohomology-comparison theorem} \label{sect_2.1.4}
From the assumption that $X$ is geometrically connected, we have
$$\mathrm{End}_\text{str}((\mathcal O_X,d)) =
\mathrm{H}_\mathrm{str}^0(X,(\mathcal O_X,d))=\bigcap_n(\mathcal O_X)^{p^n}=k.$$
On the other hand, as the underlying sheaves of objects of 
$\mathrm{str}^\mathrm{coh}(X/k)$ are locally
free, the fiber functor at any point is exact, hence is faithful.
Consequently, $\mathrm{str}^\mathrm{coh}(X/k)$ is a Tannakian category with respect to
the fiber functor at any point $x\in X(k)$. In particular, for any
stratified bundle $\mathcal{E}$, the space
$\mathrm{H}_\mathrm{str}^0(X,\mathcal{E})$ is finite dimensional. 
Tannakian duality yields an affine group scheme $\pi(X)$ (we shall not mention
the base point $x$ in what follows) with the property: the functor $x^*$ yields
an equivalence:
\begin{align}
x^*: \mathrm{str}^\mathrm{coh}(X/k)\stackrel\cong\longrightarrow \mathrm{Rep}^\mathrm{f}(\pi(X)),
\end{align}
which extends to an equivalence
\begin{align}
x^*:  \mathrm{str}^\text{ind}(X/k)\stackrel\cong\longrightarrow \mathrm{Rep}(\pi(X)).
\end{align}

According to the Tannakian duality, the functor $x^*$
yields an equivalence between $\mathrm{str}^\mathrm{ind}(X/k)$ 
and the representation category of an affine group scheme which is called the
stratified fundamental group of $X$ at the base point $x$ and is denoted
by $\pi(X,x)$:
$$x^*: \mathrm{str}^\mathrm{ind}(X/k) \cong \mathrm{Rep}(\pi(X,x)).$$
We will simplify the notation to $\pi(X)$. 
The category  $\mathrm{Rep}(\pi(X))$ has enough injectives. Hence
we can define the cohomology of $\pi(X)$ with coefficients in its representation $V$ to be
the right derived functors of the invariant-subspace functor
$$\mathrm{H}^0(\pi(X),V):=\mathrm{Hom}_{\pi(X)}
(k,V)=V^{\pi(X)},$$
where $k$ stands for the trivial representation. The $i$-th cohomology
group is denoted, as usual, by $\mathrm{H}^i(\pi(X),V)$.
We also have
$$\mathrm{H}^i(\pi(X),V)\cong \mathrm{Ext}^i_{\pi(X)}(k,V).$$
These isomorphisms yield natural maps
$$\delta^i: \mathrm{H}^i(\pi(X),E)\longrightarrow 
\mathrm{H}^i_\text{str}(X/k,\mathcal E),$$
where $E:=x^*\mathcal E=\mathcal E|_x.$ 
Since an extension in $\mathrm{str}(X/k)$ of two stratified bundles is again a stratified bundle, the maps $\delta^0$ and $\delta^1$ are isomorphisms for coherent coefficients. We shall use this degree-one comparison only in this finite-dimensional setting below.

\begin{thm}\label{thm_delta_stratified}
Let $X$ be a smooth projective curve over an algebraically closed
field $k$ of positive characteristic, $x\in X(k)$. Then for any stratified bundle 
$(\mathcal E,\nabla)$ and $E:=\mathcal E|_x$, all the maps
$$\delta^i: \mathrm{H}^i(\pi(X),E)\longrightarrow 
\mathrm{H}^i_\mathrm{str}(X/k,\mathcal E)$$
are bijective, moreover, $\delta^i$ for $i\geq 2$ are all zero-maps, that is
$$\mathrm{H}^i(\pi(X),E)=
\mathrm{H}^i_\mathrm{str}(X/k,\mathcal E)=0$$
for all $i\geq 2$. 
\end{thm}

\begin{proof}
The comparison maps in degrees zero and one are isomorphisms, as explained above. Proposition \ref{prop_finiteness} gives
\[
\mathrm H^i_{\mathrm{str}}(X/k,\mathcal E)=0\qquad(i\geq 2).
\]
It therefore suffices to prove the corresponding vanishing for $\pi(X)$.

We shall use the fact that rational group cohomology commutes with filtered direct limits (see the beginning of Section \ref{sect_comparison_ind}).

Choose an exact sequence
\[
0\longrightarrow E\longrightarrow J\longrightarrow Q\longrightarrow 0
\]
with $J$ injective in $\mathrm{Rep}(\pi(X))$. Write $Q$ as the filtered union of its finite-dimensional subrepresentations $Q_\alpha$, and let $J_\alpha\subset J$ be the inverse image of $Q_\alpha$. Then
\[
0\longrightarrow E\longrightarrow J_\alpha\longrightarrow Q_\alpha\longrightarrow 0
\]
is exact, and $\dim_k J_\alpha=\dim_k E+\dim_k Q_\alpha<\infty$.
Under Tannakian duality, it corresponds to an exact sequence of coherent stratified bundles
\[
0\longrightarrow\mathcal E\longrightarrow\mathcal J_\alpha
\longrightarrow\mathcal Q_\alpha\longrightarrow 0.
\]
Since $\mathrm H^2_{\mathrm{str}}(X/k,\mathcal E)=0$, the associated long exact sequence and the natural degree-one comparison give a surjection
\[
\mathrm H^1(\pi(X),J_\alpha)\twoheadrightarrow\mathrm H^1(\pi(X),Q_\alpha).
\]
Here the comparison is applied only to the coherent bundles $\mathcal J_\alpha$ and $\mathcal Q_\alpha$.

We have $J=\varinjlim_\alpha J_\alpha$ and $Q=\varinjlim_\alpha Q_\alpha$. Taking filtered direct limits of these surjections therefore yields
\[
0=\mathrm H^1(\pi(X),J)\twoheadrightarrow\mathrm H^1(\pi(X),Q),
\]
where the vanishing on the left follows from injectivity of $J$.
The long exact sequence of $0\to E\to J\to Q\to0$ now gives
\[
\mathrm H^2(\pi(X),E)\cong\mathrm H^1(\pi(X),Q)=0.
\]

This proves degree-two vanishing for every finite-dimensional representation. Every rational $\pi(X)$-representation is a filtered union of finite-dimensional subrepresentations, so the same commutation with filtered direct limits gives $\mathrm H^2(\pi(X),V)=0$ for every $V\in\mathrm{Rep}(\pi(X))$.
For such a $V$, choose $0\to V\to I\to W\to0$ with $I$ injective. Dimension shifting gives
\[
\mathrm H^i(\pi(X),V)\cong\mathrm H^{i-1}(\pi(X),W)\qquad(i\geq 2).
\]
Induction on $i$, starting with degree two and applying the assertion to every rational representation, proves $\mathrm H^i(\pi(X),V)=0$ for all $i\geq 2$. Taking $V=E$ completes the proof.
\end{proof}

\subsection{Comparison with cohomology of the \'etale fundamental group}
\label{sect_etale_trivial}
We call a vector bundle on $X$  \textit{\'etale finite}
if it can be trivialized by an \'etale finite torsor 
(i.e. torsor under an \'etale finite group scheme) on $X$.
The category $\mathcal C^\mathrm{\'et}(X)$
of \'etale finite vector bundles equipped with the fiber functor $x^*$ is Tannakian (the fixed $k-$ point $x$ will be omitted) and its dual
is precisely the (geometric) \'etale fundamental group scheme $\pi^\mathrm{\'et}(X)$ (considered as
a pro-finite group scheme).
The reader is referred to \cite{No82}, \cite{MS02} and \cite{EHS08} for more details. 
Further, there exists a natural functor 
$$\mathcal C^\mathrm{\'et}(X)\longrightarrow \mathrm{str}(X/k)$$
in the sense that each \'etale finite vector bundle is equipped in a canonical way
a stratification making it a stratified bundle with finite monodromy (cf. \cite{dSa07}). This
functor is fully faithful, exact and closed under taking subquotients (which is
equivalent to the surjectivity of the map $\pi(X)\longrightarrow \pi^\mathrm{\'et}(X)$).
For convenience, we shall identify $\mathcal C^\mathrm{\'et}(X)$ with
the full subcategory of stratified bundle with finite monodromy
$\mathrm{str}^\text{fin}(X/k)$. 

We notice that the condition ``\'etale finite'' is equivalent to the condition ``\'etale trivializable'' as 
defined in \cite{LS77}. We shall use interchangeably these notations.

Let $\mathrm{str}^\mathrm{uni}(X)$ be the full subcategory of 
$\mathrm{str}^\mathrm{coh}(X)$ of objects which are iterated extension of the
trivial object $(\mathcal O,d)$. 
The corresponding Tannakian group is a pro-unipotent group, denoted $\pi^{\rm uni}(X).$
It is a quotient group scheme of $\pi(X)$, and, according to \cite{dSa07}, it is
pro-finite, hence is a quotient of $\pi^\mathrm{\'et}(X).$

Summarizing, we have surjective (i.e. faithfully flat) maps
$$\pi(X)\longrightarrow \pi^\mathrm{\'et}(X)\longrightarrow
\pi^\mathrm{uni}(X).$$

\begin{prop}\label{prop_compare_etale} Let $X/k$ be a smooth projective curve.
For any finite representation $V$ of $\pi^\mathrm{\'et}(X)$ the natural maps induced by restriction
$$\mathrm{res}^i:\mathrm H^i(\pi^\mathrm{\'et}(X),V)\longrightarrow 
\mathrm H^i(\pi(X),V)$$
are isomorphisms. 
\end{prop}
\begin{proof}
The claim for $i=0$ follows from the surjectivity of $\pi(X)\longrightarrow \pi^\mathrm{\'et}(X)$. The surjectivity also ensures the injectivity of $\mathrm{res}^1$.

The claim for $i=1$ amounts to the thickness of 
$\mathcal C^\mathrm{\'et}(X)$ in $ \mathrm{str}(X/k)$, that is, for any extension
$$0 \longrightarrow \mathcal V\longrightarrow \mathcal W\longrightarrow \mathcal O_X \longrightarrow 0$$
in $ \mathrm{str}(X/k)$ with $\mathcal V$ in $\mathrm{str}^\text{fin}(X/k)$ we also
have $\mathcal W$ in $\mathrm{str}^\text{fin}(X/k)$.
Indeed, let $Y\stackrel p\longrightarrow X$ be the $G$-torsor constructed from $\mathcal V$.
By assumption, $G$ is an \'etale finite group scheme. Then $p^*\mathcal V$ is 
a trivial vector bundle on $Y$ equipped with the trivial connection. 
According to dos Santos' result mentioned above ($\pi^{\rm uni}(X)$ is a quotient of $\pi^\mathrm{\'et}(X)$), $p^*\mathcal W$ is also a
finite connection on $Y$. But this implies the push forward $p_*p^*\mathcal W$
is also a finite connection that forces
$\mathcal W$ to be a finite connection as we have an inclusion of connections
$$\mathcal W\hookrightarrow  p_*p^*\mathcal W.$$ Thus $\mathrm{res}^1$ is an
isomorphism. 

For $i\geq 2,$ we use the same trick as in the proof of \ref{thm_delta_stratified}.
First notice that the isomorphism $\mathrm{res}^1$ extends to any representations 
of $\pi^\mathrm{\'et}(X)$.
Let $J$ be the injective envelope of $V$ in $\mathrm{Rep}(\pi^\mathrm{\'et}(X))\cong
\mathrm{Ind}-\mathcal C^\mathrm{\'et}(X)$. Then we have
$$\mathrm{H}^2(\pi^\mathrm{\'et}(X),V)\cong 
\mathrm{H}^1(\pi^\mathrm{\'et}(X),J/V)\stackrel{(\mathrm{res}^1)}\cong 
\mathrm{H}^1(\pi(X),J/V)\cong
\mathrm{H}^2(\pi(X),V)=0.$$
The proof is complete. 
\end{proof}

As a corollary we obtain the old result due to Serre and Pacheco-Stevenson \cite{Ser90,PS00}. 
\begin{cor}The cohomology groups 
$\mathrm H^i(\pi^\mathrm{\'et}(X),V)$ vanish for $i\geq 2$.\end{cor}
\subsection{Freeness of the pro-unipotent fundamental group scheme}
By definition, the cohomology of $\pi^{\rm uni}(X)$ with coefficients in a representation
$V$ is the same as the cohomology of $(\mathcal V,\nabla)$ in  $\mathrm{str}^\mathrm{uni}(X/k)$
with respect to its ind-category, where $(\mathcal V,\nabla)$ is the nilpotent
stratified bundle corresponding to the representation $V$:
$$\mathrm H^i(\pi^\mathrm{uni}(X),V)\cong \mathrm H^i_{\mathrm{str}^\mathrm{uni}}(X/k,(\mathcal V,\nabla)).$$
Using Tannakian duality and the fact that extensions of nilpotent stratified bundles
in $\mathrm{str}(X/k)$ are again nilpotent, we have, for $i=0,1$,
$$\mathrm H^i_{\mathrm{str}^\mathrm{uni}}(X/k,(\mathcal V,\nabla))\cong \text H^i_\mathrm{str}(X/k,(\mathcal V,\nabla))$$
and 
$$\mathrm H^i(\pi^\mathrm{uni}(X),V)\cong\mathrm  H^i(\pi(X),V)$$
for any nilpotent stratified bundle $(\mathcal V,\nabla)$ and its fiber $V$. 
Recall that the latter isomorphism extends to ind-representations, while the former
one may not extend to ind-stratified bundles. 
\begin{cor}\label{cor_nilpotent_cohomology}
The group cohomology 
$\mathrm  H^i(\pi^\mathrm{uni}(X),V)$, $i\geq 2$, vanishes for all representations $V$.
\end{cor} 
\begin{proof} 
Let $J$ be an injective envelope of $V$ in $\mathrm{Rep}(\pi^\mathrm{uni}(X/k))$.
Then we have a natural isomorphism  (cf. proof of Proposition \ref{prop_compare_etale}):
$$\mathrm H^2(\pi^\mathrm{uni}(X),  V)\cong 
\mathrm H^{1}(\pi^\mathrm{uni}(X),J/V).$$
Employing the isomorphism
$$\mathrm H^1(\pi^\mathrm{uni}(X),V)\cong\mathrm  H^1(\pi(X),V)$$
which holds for any representation of $\pi^\mathrm{uni}(X)$, we obtain an isomorphism
$$\mathrm H^2(\pi^\mathrm{uni}(X),  V)\hookrightarrow \mathrm H^2(\pi(X),  V).$$
Whence we obtain the vanishing claim, by means of Theorem \ref{thm_delta_stratified}. The case of $i\geq 3$ follows by a similar argument.
\end{proof}
We now introduce the notion of freeness in the category of pro-\'etale
pro-unipotent group schemes. This is the natural analogue, on the group-scheme
side, of the usual universal property of a free pro-\(p\) group.

\begin{defn}
A pro-\'etale pro-unipotent affine group scheme over \(k\) is an affine group
	scheme \(G\) which can be written as a projective limit
	\[
	G=\varprojlim_i G_i,
	\]
	where each \(G_i\) is a finite \'{e}tale unipotent group scheme over \(k\).
	Equivalently, since \(k\) is algebraically closed, the group schemes \(G_i\)
	are finite constant \(p\)-group schemes.
\end{defn}

\begin{defn} \label{free pro}
	Let \(S\) be a set. A pro-\'etale pro-unipotent affine group scheme
	\(F^{\mathrm{uni,et}}_S\) over \(k\) is called the free pro-\'etale
	pro-unipotent group scheme generated by \(S\) if, for every pro-\'etale
	pro-unipotent affine group scheme \(G\) over \(k\), there is a natural
	bijection
	\[
	\operatorname{Hom}_{k\text{-grp.sch}}
	\bigl(F^{\mathrm{uni,et}}_S,G\bigr)
	\simeq
	\operatorname{Map}(S,G(k)).
	\]
\end{defn}
We first establish the equivalence between the category of pro-\(p\) groups
and the category of pro-\'etale pro-unipotent affine group schemes over \(k\).
\begin{prop}\label{prop:proetale-unipotent-prop-groups}
	The
	functor
	\[
	P \longmapsto \underline{P}_k
	\]
	from pro-\(p\) groups to pro-\'{e}tale pro-unipotent affine group schemes over
	\(k\) is an equivalence of categories. Its quasi-inverse is given by
	\[
	G \longmapsto G(k).
	\]
\end{prop}

\begin{proof}
	We first prove the assertion at the finite level. Since \(k\) is algebraically closed, every finite \'{e}tale \(k\)-scheme is a
	finite disjoint union of copies of \(\operatorname{Spec} k\)
	\cite[Tag 04GL]{Stack}. Therefore every finite \'{e}tale group scheme
	\(H\) over \(k\) is constant:
	\[
	H\simeq \underline{H(k)}_k .
	\]
	Thus the category of finite \'{e}tale group schemes over \(k\) is equivalent to
	the category of finite groups.
	
	Under this equivalence, the representation category of
	\(\underline{\Gamma}_k\) is the category of \(k[\Gamma]\)-modules. Hence
	\(\underline{\Gamma}_k\) is unipotent if and only if every simple
	\(k[\Gamma]\)-module is trivial. If \(\Gamma\) is a finite \(p\)-group, then
	\(k[\Gamma]\) is a local algebra; equivalently, it has a unique simple module,
	namely the trivial one. Conversely, if \(\Gamma\) is not a \(p\)-group, choose a
	prime \(\ell\neq p\) dividing \(|\Gamma|\). By Cauchy's theorem, \(\Gamma\)
	contains a cyclic subgroup \(C\simeq \mathbb Z/\ell\mathbb Z\). Since \(k\) is
	algebraically closed and \(\ell\neq p\), \(C\) has a non-trivial
	one-dimensional character over \(k\). Inducing this character to \(\Gamma\) and
	taking a simple quotient gives a non-trivial simple \(k[\Gamma]\)-module. Hence
	\(\underline{\Gamma}_k\) is unipotent if and only if \(\Gamma\) is a finite
	\(p\)-group.
	
	Now let \(P=\varprojlim_i P_i\) be a pro-\(p\) group, with \(P_i\) finite
	\(p\)-groups. Define
	\[
	\underline{P}_k:=\varprojlim_i \underline{P_i}_k .
	\]
	This is a pro-\'{e}tale pro-unipotent affine group scheme. Conversely, if
	\(G=\varprojlim_i G_i\) is pro-\'{e}tale pro-unipotent, then each \(G_i\) is
	finite \'{e}tale unipotent, hence
	\[
	G_i\simeq \underline{G_i(k)}_k
	\]
	with \(G_i(k)\) a finite \(p\)-group. Therefore
	\[
	G(k)=\varprojlim_i G_i(k)
	\]
	is a pro-\(p\) group and
	\[
	G\simeq \varprojlim_i \underline{G_i(k)}_k
	= \underline{G(k)}_k .
	\]
	It remains to prove that the functor is fully faithful. Write
	\[
	P=\varprojlim_i P_i,\qquad Q=\varprojlim_j Q_j
	\]
	as inverse limits of finite \(p\)-groups. Since
	\(\underline{Q}_k=\varprojlim_j \underline{Q_j}_k\), a morphism
	\(\underline{P}_k\to \underline{Q}_k\) is the same thing as a compatible family
	of morphisms
	\[
	\underline{P}_k\longrightarrow \underline{Q_j}_k .
	\]
	For fixed \(j\), the group scheme \(\underline{Q_j}_k\) is finite \'{e}tale,
	hence of finite presentation over \(k\). Since
	\(\underline{P}_k=\varprojlim_i \underline{P_i}_k\) is an inverse limit with
	affine transition morphisms, every morphism
	\(\underline{P}_k\to \underline{Q_j}_k\) factors through some finite stage
	\(\underline{P_i}_k\), by the standard limit argument for morphisms into a
	finitely presented scheme \cite[Tag 01ZM]{Stack}. Hence
	\[
	\begin{aligned}
		\operatorname{Hom}_{k\text{-grp.sch}}
		\bigl(\underline{P}_k,\underline{Q}_k\bigr)
		&\simeq
		\varprojlim_j
		\operatorname{Hom}_{k\text{-grp.sch}}
		\bigl(\underline{P}_k,\underline{Q_j}_k\bigr)     \\
		&\simeq
		\varprojlim_j \varinjlim_i
		\operatorname{Hom}_{k\text{-grp.sch}}
		\bigl(\underline{P_i}_k,\underline{Q_j}_k\bigr)   \\
		&\simeq
		\varprojlim_j \varinjlim_i
		\operatorname{Hom}_{\mathrm{grp}}(P_i,Q_j).
	\end{aligned}
	\]
	On the other hand, a continuous homomorphism \(P\to Q\) is the same thing as a
	compatible family of continuous homomorphisms \(P\to Q_j\). Since \(Q_j\) is
	finite and discrete, each such homomorphism factors through a finite quotient
	\(P_i\). Thus
	\[
	\operatorname{Hom}_{\mathrm{cont.grp}}(P,Q)
	\simeq
	\varprojlim_j \varinjlim_i
	\operatorname{Hom}_{\mathrm{grp}}(P_i,Q_j).
	\]
	This proves full faithfulness, and hence the proposition.
\end{proof}

\begin{lem}\label{lem:cohomology-comparison-proconstant}
	Let \(P\) be a pro-\(p\) group and let
	\(\underline{P}_k\) be the associated pro-constant affine group scheme over
	\(k\). For every discrete continuous \(k\)-linear representation \(V\) of
	\(P\), viewed equivalently as a representation of \(\underline{P}_k\), there is
	a natural isomorphism
	\[
	\mathrm H^n(\underline{P}_k,V)
	\simeq
	\mathrm H^n_{\mathrm{cts}}(P,V).
	\]
	In particular,
	\[
	\mathrm H^n(\underline{P}_k,k)
	\simeq
 \mathrm H^n_{\mathrm{cts}}(P,k).
	\]
	Moreover, if \(M\) is a finite discrete \(\mathbb F_p[P]\)-module, then
	\[
	\mathrm H^n_{\mathrm{cts}}(P,M)\otimes_{\mathbb F_p}k
	\simeq
	\mathrm H^n_{\mathrm{cts}}\bigl(P,M\otimes_{\mathbb F_p}k\bigr).
	\]
\end{lem}

\begin{proof}
	The category \(\operatorname{Rep}_k(P_k)\) identifies with the category of
	discrete continuous \(k\)-linear representations of \(P\). Under this
	identification, the invariant functor is the usual invariant functor
	\[
	V\longmapsto V^P .
	\]
	Therefore its right derived functors agree with continuous group cohomology:
	\[
	\mathrm H^n(P_k,V)=R^n(V\mapsto V^P)= \mathrm H^n_{\mathrm{cts}}(P,V).
	\]
	Equivalently,
	\[
	\mathrm H^n(P_k,V)
	=
	\operatorname{Ext}^n_{\operatorname{Rep}_k(P_k)}(k,V)
	\simeq
	\mathrm H^n_{\mathrm{cts}}(P,V).
	\]
	
	For the last assertion, compute \(\mathrm H^n_{\mathrm{cts}}(P,M)\) by the standard
	complex of continuous inhomogeneous cochains
	\[
	C^q_{\mathrm{cts}}(P,M)=\operatorname{Cont}(P^q,M).
	\]
	Since \(M\) is finite and \(k\) is flat over \(\mathbb F_p\), tensoring the
	continuous cochain complex with \(k\) gives the cochain complex computing
	\[
	\mathrm H^\bullet_{\mathrm{cts}}(P,M\otimes_{\mathbb F_p}k).
	\]
	Hence
	\[
	\mathrm H^n_{\mathrm{cts}}(P,M)\otimes_{\mathbb F_p}k
	\simeq
	\mathrm H^n_{\mathrm{cts}}(P,M\otimes_{\mathbb F_p}k).
	\]
\end{proof}

\begin{prop}\label{thm:pi-uni-free}
	The group scheme \(\pi^{\mathrm{uni}}(X)\) is a free pro-\'etale
	pro-unipotent group scheme on
	$
	h^1(\pi^{\mathrm{uni}},k)
	$
	generators.
\end{prop}

\begin{proof}
	We apply Proposition~\ref{prop:proetale-unipotent-prop-groups} to
	$
	\pi^{\mathrm{uni}}:=\pi^{\mathrm{uni}}(X)
	$ which is a pro-\'etale pro-unipotent group scheme (see Section \ref{sect_etale_trivial}). Hence there
	exists a pro-\(p\) group \(P\) such that
	\[
	\pi^{\mathrm{uni}}\simeq \underline{P}_k.
	\]
	
	By Lemma~\ref{lem:cohomology-comparison-proconstant}, the vanishing
	\[
	\mathrm H^2(\pi^{\mathrm{uni}},W)=0
	\]
	for every rational representation \(W\) of \(\pi^{\mathrm{uni}}\) implies, in
	particular, that
	\[
	\mathrm H^2_{\mathrm{cts}}(P,M\otimes_{\mathbb F_p}k)=0
	\]
	for every finite discrete \(\mathbb F_p[P]\)-module \(M\). By the base-change
	isomorphism in Lemma~\ref{lem:cohomology-comparison-proconstant}, this gives
	\[
	\mathrm H^2_{\mathrm{cts}}(P,M)\otimes_{\mathbb F_p}k=0.
	\]
	Since \(k/\mathbb F_p\) is faithfully flat, it follows that
	\[
	\mathrm H^2_{\mathrm{cts}}(P,M)=0
	\]
	for every finite discrete \(\mathbb F_p[P]\)-module \(M\). Hence
	\[
	\operatorname{cd}_p(P)\leq 1.
	\]
	By Serre's criterion for pro-\(p\) groups
	\cite[Chapter I, Section 4]{SerreGalois}, a pro-\(p\) group has
	\(p\)-cohomological dimension at most \(1\) if and only if it is free
	pro-\(p\). Therefore \(P\) is free pro-\(p\).
	
	Let
	\[
	r:=\dim_k \mathrm H^1(\pi^{\mathrm{uni}},k).
	\]
	Then
	\[
	\begin{aligned}
		r
		&=
		\dim_k \mathrm H^1(\pi^{\mathrm{uni}},k) \\
		&=
		\dim_k\bigl(\mathrm H^1_{\mathrm{cts}}(P,\mathbb F_p)\otimes_{\mathbb F_p}k\bigr)\\
		&=
		\dim_{\mathbb F_p}\mathrm H^1_{\mathrm{cts}}(P,\mathbb F_p).
	\end{aligned}
	\]
	For a free pro-\(p\) group, this number is precisely the minimal number of
	topological generators. Hence
	\[
	P\simeq F^{(p)}_r,
	\]
	where \(F^{(p)}_r\) denotes the free pro-\(p\) group on \(r\) generators.
	Consequently
	\[
	\pi^{\mathrm{uni}}\simeq \underline{F^{(p)}_r}_k .
	\]
	It remains to express this conclusion intrinsically in the language of group
	schemes (see Definition \ref{free pro}). Let \(S\) be a set with \(|S|=r\). For every pro-\'etale
	pro-unipotent affine group scheme \(G\) over \(k\), write \(G\simeq \underline{Q}_k\), where
	\(Q=G(k)\). Since \(P\mapsto P_k\) is an equivalence of categories, we get
	\[
	\begin{aligned}
		\operatorname{Hom}_{k\text{-grp.sch}}(\pi^{\mathrm{uni}},G)
		&\simeq
		\operatorname{Hom}_{k\text{-grp.sch}}\bigl( \underline{F^{(p)}_r}_k,\underline{Q}_k\bigr) \\
		&\simeq
		\operatorname{Hom}_{\mathrm{cont.grp}}\bigl(F^{(p)}_r,Q\bigr) \\
		&\simeq
		\operatorname{Map}(S,Q) \\
		&=
		\operatorname{Map}(S,G(k)).
	\end{aligned}
	\]
	Thus \(\pi^{\mathrm{uni}}\) is the free pro-\'etale pro-unipotent group scheme
	on \(r=h^1(\pi^{\mathrm{uni}},k)\) generators.
\end{proof}
\subsection{Application: a computation for cohomology of curves of genus $1$} 
Let $X$ have genus $1$. 
Then, according to dos Santos
\cite{dSa07}, $\pi(X)$ is commutative and decomposes into a direct product of a unipotent part and a diagonal part. 
 Since the group cohomology of the diagonal part
vanishes (being linearly reductive), 
$$\mathrm{H}^i(\pi(X),V)=
\mathrm{H}^i(\pi^\mathrm{uni}(X),V),$$
 for $i\geq 1$ and 
any unipotent representation $V$. In categorical terms, we have: 
\begin{itemize}
\item If $X$ is supersingular, the $p$-rank is zero, so the unipotent part of $\pi(X)$ is trivial. Hence $\pi(X)$ is diagonal.
\item If $X$ is ordinary, that is, the absolute Frobenius acts by a bijective map, 
then according to \cite{dSa07}, any object of $\mathrm{str}^\mathrm{coh}(X/k)$ decomposes
into a direct sum of objects of the form
$$\mathcal L\otimes\mathcal  N$$
where $\mathcal L$ is a rank one stratified bundle and $\mathcal N$ is a nilpotent
object, i.e., all of its irreducible subquotients
are trivial objects. Such an object is also called $\mathcal L$-isotypical.
\end{itemize}

Now, dos Santos tells us that (loc.cit. Theorem 21)
$$\pi^\mathrm{uni}(X)=\mathbb Z_p^r,$$
where $r$ is the $p$-rank of $X$ (i.e., the rank of the Frobenius map acting on 
$\mathrm{H}^1(X,\mathcal O_X)$, which can be either $1$ or $0$). 
Thus we have a complete understanding of the stratified cohomology of stratified bundles on an elliptic curve.
\begin{prop}\label{prop_vanishing_elliptic}
Let $X$ be a smooth projective curve of genus 1. 
\begin{enumerate}
\item If $X$ has $p$-rank 0 then $\mathrm{h}^1_\mathrm{str}(X,\mathcal E)=0$ for any
stratified bundle $\mathcal E$. 
\item If $X$ has $p$-rank 1, then $\mathrm{h}^1_\mathrm{str}(X,\mathcal O_X)=1$ and hence this function is non-zero at any unipotent stratified bundle. On the other hand,
for any non-trivial stratified line bundle $\mathcal L$, we have $\mathrm{h}^1_\mathrm{str}(X,\mathcal V)=0$ if $\mathcal V$ is an $\mathcal L$-isotypical stratified bundle.
\end{enumerate}
\end{prop}
\begin{proof}
If $X$ has $p$-rank 0, that is, $X$ is supersingular, the action of the Frobenius of $X$ is zero,
hence  $\pi(X)$ is diagonal. So its cohomology vanishes from degree 1. 
 
If $X$ has $p$-rank 1, that is $X$ is ordinary, the absolute Frobenius acts by a bijective map
on the 1-dimensional space $\mathrm{H}^1(X,\mathcal O_X)$. 
Hence $\mathrm{H}^1_\mathrm{str}(X,\mathcal O_X)\cong k$. If $\mathcal V$ is a nilpotent object, then there is a surjective map $\mathcal V\longrightarrow \mathcal O_X$ as objects in $\mathrm{str}(X/k)$. As the $\mathrm{H}^2_\mathrm{str}$-group vanishes, we have
a surjective map $\mathrm{H}^1_\mathrm{str}(X,\mathcal V)\longrightarrow 
\mathrm{H}^1_\mathrm{str}(X,\mathcal O_X)$, whence $\mathrm{H}^1_\mathrm{str}(X,\mathcal V)\neq 0$.

On the other hand, as explained above, any object of $\mathrm{str}^\mathrm{coh}(X/k)$ decomposes
into a direct sum of objects of the form
$\mathcal L\otimes\mathcal  N$
where $\mathcal L$ is a rank one stratified bundle and $\mathcal N$ is a nilpotent
object, i.e., all of its irreducible subquotients
are trivial. Hence, if  $\mathcal L$ is a non-trivial rank 1 object
$\mathrm{H}^1_\mathrm{str}(X,\mathcal L)\cong 
\mathrm{Ext}^1_{\mathrm{str}(X/k)}(\mathcal O_X,\mathcal L)=0.$
Using the long exact sequence of cohomology groups, we obtain 
$\mathrm{Ext}^1_{\mathrm{str}(X/k)}(\mathcal N,\mathcal L)=0$ for any nilpotent stratified bundle
$\mathcal N$, whence the last claim of Proposition. 
\end{proof}
\begin{cor}
Consider $\mathbb Z_p$ as a pro-finite group scheme over $k$. Then 
$\mathrm{h}^i(\mathbb Z_p, -)$ vanishes for any $i\geq 2$ and
$$\mathrm{h}^i(\mathbb Z_p, k)=1,$$
for $i=0,1$.  
\end{cor}

\subsection{The non-vanishing of the first cohomology group}
\label{sect_comparison_torsion}

Our aim in this section is to discuss the vanishing property of the first cohomology group.  
This question is motivated by the non-vanishing property for the first de Rham
cohomology on curves over a field of characteristic zero \cite{BHT25}.
We shall restrict ourselves to \textit{finite stratified bundles}, that is
\'etale finite (i.e. trivializable) vector bundles, cf. \ref{sect_etale_trivial}.

A vector bundle is said to be $F$-periodic if it is isomorphic to the pull-back of itself by
some powers of the Frobenius map. 
According to Lange-Stuhler \cite{LS77}, $F$-periodic vector bundles are \'etale trivializable.
The converse implication holds if the base field is finite or if the
vector bundle is stable \cite{BD07}. Now we restrict further our problem to $F$-periodic vector bundles.  

Let $\mathcal E$ be an $F$-periodic vector bundle, we shall use the same symbol for
the associated stratified bundle. Fix an isomorphism 
$\tau: F^{*n}\mathcal E\longrightarrow
\mathcal E$.  
The composition of $\tau$'s and the canonical inclusion
$\mathcal E\longrightarrow F^{*n}  \mathcal E $ yields a $p^n$-linear endomorphism
$$f:\mathcal E \longrightarrow  F^{*n}  \mathcal E\xrightarrow{\,\,\,\tau\,\,\,} \mathcal E.$$ 
Denote by $\mathcal H_{\mathcal E}$ the subsheaf of $\mathcal E$ 
fixed by $f$ in the \'etale topology on $X$. The lemma below goes back to
 \cite[Propositions 1.1 and 1.2]{Kat69}.
 
\begin{lem}\label{lem_etale_sheaf}
Let \(\mathcal E\) be an \(F\)-periodic stratified bundle of period \(n\).
Then $\mathcal H_{\mathcal E}$ is a locally constant sheaf with fiber a vector space over $\mathbb F_{p^n}$ with dimension equal to the rank of $\mathcal E$, 
and there is an exact sequence 
$$0 \rightarrow \mathcal H_{\mathcal E} \longrightarrow 
\mathcal E  \xrightarrow{f-1} \mathcal E  \longrightarrow 0$$ 
of sheaves on $X_\text{\'et}.$
\end{lem}

\subsubsection*{Remark} $\mathcal H_{\mathcal E}$ can be seen as a finite additive group scheme on $X$.

\begin{prop} \label{etale comparison}
Let \(\mathcal E\) be an \(F\)-periodic stratified bundle of period \(n\), and let
\(\mathcal H_{\mathcal E}\) be the associated finite \'{e}tale group scheme defined above. Then we have a canonical isomorphism:
$$\mathrm H^i_\text{\'et}(X, \mathcal H_{\mathcal E}) \otimes_{\mathbb{F}_q} k \xrightarrow{\cong} \mathrm H^i_\mathrm{str}(X, \mathcal E),$$ for all $i,$ where $q = p^n.$
\end{prop}
\begin{proof}
From the short exact sequence of sheaves on $X_\text{\'et}:$ 
$$0 \rightarrow \mathcal H_{\mathcal E} \rightarrow\mathcal E  \xrightarrow{F^{*n}-1} 
\mathcal E  \rightarrow 0,$$
we have the following long exact sequence of \textit{\'etale} cohomologies:
\begin{equation}
\cdots \rightarrow\mathrm  H^i_\text{\'et}(X, \mathcal H_{\mathcal E}) 
\rightarrow\mathrm  H^i_\text{\'et}(X,\mathcal E) \xrightarrow{F^{*n}-1}
\mathrm  H^i_\text{\'et}(X, \mathcal E) \rightarrow \cdots
\end{equation}
Since $\mathcal E$ is a coherent sheaf on $X,$ we have 
$\mathrm H^i_\text{\'et}(X, \mathcal E) \cong \mathrm H^i(X, \mathcal E).$ 
Since the induced morphism on cohomology 
$F^n:\mathrm H^i_\text{\'et}(X, \mathcal E)  \rightarrow 
\mathrm  H^i_\text{\'et}(X, \mathcal E)$ 
is $q$-linear, the argument before this theorem implies that 
$F^{*n}-1: \mathrm  H^i_\text{\'et}(X,\mathcal  E)  \rightarrow 
\mathrm  H^i_\text{\'et}(X,\mathcal  E)$ is surjective
 (cf. \cite[Proposition~1.2]{Kat69}). 
Hence, for each $i,$ we have the following short exact sequence:
$$ 0 \rightarrow \mathrm  H^i_\text{\'et}(X,\mathcal H_{\mathcal E}) \rightarrow 
\mathrm H^i_\text{\'et}(X,\mathcal E)  \xrightarrow{F^{*n}-1} 
\mathrm  H^i_\text{\'et}(X,\mathcal E) \rightarrow 0.$$
This finishes the proof since we have 
$\mathrm H^i_\mathrm{str}(X, \mathcal E) =\mathrm  H^i(X, \mathcal E)_\mathrm{ss} \cong
\mathrm  H^i(X, \mathcal E)^{F^{*n}-1}
 \otimes_{\mathbb F_q}k $ which is isomorphic to 
$\mathrm  H^i_\text{\'et}(X, \mathcal H_{\mathcal E}) \otimes_{\mathbb F_q}k$ by the above short exact sequence.
\end{proof}

Assume that $X$ is a curve of genus larger than one. The method of \cite{BHT25} relies
on the Euler characteristic formula which is not always available in positive characteristic.
More precisely, let $\mathcal F$ be a locally constant sheaf on $X$ and 
$f:Y\longrightarrow X$
be an \'etale covering map that trivializes $\mathcal F$, that is, $f^*\mathcal F$ is
a constant sheaf on $Y$. Then it holds \cite[Theorem~0.2]{Pin00}:
$$\chi(X,\mathcal F)\geq (1-g_X)\cdot \mathrm{rank\ }\mathcal F,$$
where $\mathrm{rank\ }\mathcal F$ denotes the ${\mathbb F}_{p}$-dimension 
of the fiber of $\mathcal F$ at a point of $X$.  
The equality takes place when $Y$ is ordinary. We introduce the following notation.
\begin{defn}
A vector bundle $\mathcal{F}$ on a smooth projective curve
$X$ is said to be ordinarily \'etale trivializable if there 
exists a Galois \'etale covering map $f:Y\longrightarrow X$
that trivializes $\mathcal F$, where $Y$ is an ordinary curve.
\end{defn}

\begin{cor}\label{cor_non-vanishing}
Let $X$ be an ordinary curve of genus larger than $1$. 
Let $\mathcal E$ be an ordinarily \'etale trivializable vector bundle. 
Then $\mathrm h^1_\mathrm{str}(X, \mathcal E) \neq 0.$
\end{cor}
\begin{proof}
\'Etale trivializable vector bundles are semi-stable of degree 0, hence possess
composition series  such that the successive quotients are 
stable of degree zero.
Using the long exact sequence of cohomology and induction on the length
of the above composition series, the problem is essentially reduced
to considering only \'etale trivializable \textit{stable} vector bundles.
According to  \cite{BD07}, \'etale trivializable stable vector bundles are $F$-periodic.
Thus, we can assume that $\mathcal E$ is $F$-periodic with $F$-period, say, $n$.

By Proposition \ref{etale comparison}, it is enough to show that 
$\mathrm h^1_\text{\'et}(X,\mathcal H_{\mathcal E})\neq 0$.  By \cite[Theorem~0.2]{Pin00}, we have that:
\begin{align*}
\chi(H_\mathcal{E})= (1-\mathrm h^1_\mathrm{str}(X,\mathcal{O}_X))\cdot \mathrm{rank}
\mathcal H_{\mathcal E}\neq 0.
\end{align*} This proves our claim.
\end{proof}

We shall give in the Appendix \ref{proof_non-example} an
example of an $F$-periodic vector bundle whose first stratified cohomology vanishes. 

\section{Comparison for cohomologies of ind-stratified bundles}
\label{sect_comparison_ind}
For any affine group scheme $G$ over $k$, rational group cohomology commutes with filtered direct limits:
\[
\mathrm H^i\bigl(G,\varinjlim_\alpha V_\alpha\bigr)
\cong\varinjlim_\alpha\mathrm H^i(G,V_\alpha).
\]
Indeed, the standard Hochschild complex computing this cohomology has terms
$C^n(G,V)=V\otimes_k k[G]^{\otimes n}$ (cf. \cite[Part~I, Chapter~4]{Jan87}).
These commute with filtered direct limits, which are exact over $k$, so taking cohomology commutes with these limits as well.
In contrast, stratified cohomology involves inverse limits and need not commute with filtered direct limits. Consequently, the comparison maps of Theorem \ref{thm_delta_stratified} need not be isomorphisms for ind-stratified bundles, as the examples below show.

\subsection{The case of elliptic curves}
We first discuss the comparison problem for ind-stratified bundles on elliptic curves. The answer depends on whether the curve is supersingular or ordinary.

\begin{prop}\label{prop_elliptic}
Let $X/k$ be an elliptic curve.
\begin{enumerate}
\item If $X$ is supersingular, then the maps $\delta^i$ defined as in Theorem \ref{thm_delta_stratified} are isomorphisms for every ind-stratified bundle and every $i\geq 0$.
\item If $X$ is ordinary, then the same conclusion holds for every $\mathcal L$-isotypical ind-stratified bundle, where $\mathcal L$ is a rank-one stratified bundle. Consequently, it also holds for finite direct sums of isotypical ind-stratified bundles.
\end{enumerate}
\end{prop}

\begin{proof}
Assume first that $X$ is supersingular. By \cite[Theorem~21]{dSa07}, the unipotent part of $\pi(X)$ is trivial, hence $\pi(X)$ is diagonal. Thus every ind-stratified bundle is a direct sum of rank-one objects.

Let $\mathcal L$ be a non-trivial stratified line bundle. Since $\operatorname{Pic}^0(X)[p](k)=0$, multiplication by $p$ is injective on $\operatorname{Pic}^0(X)(k)$. It follows that $\mathcal L^{(n)}$ is non-trivial for every $n$: indeed, if one term were trivial, then all terms would be trivial and hence so would $\mathcal L$. Therefore
\[
\mathrm H^0(X,\mathcal L^{(n)})=\mathrm H^1(X,\mathcal L^{(n)})=0
\qquad\text{for all }n.
\]
For the trivial character, the transition maps on $\mathrm H^0(X,\mathcal O_X)$ are bijective, whereas those on $\mathrm H^1(X,\mathcal O_X)$ are zero. Ogus' exact sequence \cite[Theorem~2.4]{Og75} therefore gives
\[
\mathrm H^i_{\mathrm{str}}(X/k,\mathcal V)=0\qquad(i\geq1)
\]
for every ind-stratified bundle $\mathcal V$. The same vanishing holds for group cohomology because $\pi(X)$ is diagonal.

Now assume that $X$ is ordinary and let $\mathcal V=\mathcal L\otimes\mathcal N$ be $\mathcal L$-isotypical, with $\mathcal N$ ind-unipotent. If $\mathcal L$ is trivial, Frobenius acts bijectively on $\mathrm H^i(X,\mathcal O_X)$ for $i=0,1$. By induction on the length of a finite unipotent object, and then by passage to filtered direct limits, all transition maps
\[
\mathrm H^i(X,\mathcal N^{(n+1)})\longrightarrow
\mathrm H^i(X,\mathcal N^{(n)}),\qquad i=0,1,
\]
are bijective. Hence the first derived inverse limits vanish and
\[
\mathrm H^1_{\mathrm{str}}(X/k,\mathcal N)
\cong \mathrm H^1(X,\mathcal N)
\cong \varinjlim_\alpha \mathrm H^1_{\mathrm{str}}(X/k,\mathcal N_\alpha),
\]
where $(\mathcal N_\alpha)$ is the system of coherent unipotent subobjects. Since group cohomology commutes with filtered direct limits, Theorem \ref{thm_delta_stratified} gives the comparison in this case.

If $\mathcal L$ is non-trivial, there exists $n_0$ such that $\mathcal L^{(n)}$ is a non-trivial degree-zero line bundle for every $n\geq n_0$. A non-trivial degree-zero line bundle on an elliptic curve has vanishing cohomology. Using finite unipotent filtrations and then filtered direct limits, we obtain
\[
\mathrm H^i\bigl(X,\mathcal L^{(n)}\otimes\mathcal N^{(n)}\bigr)=0,
\qquad i=0,1,\quad n\geq n_0.
\]
Thus the inverse systems occurring in Ogus' exact sequence \cite[Theorem~2.4]{Og75} are eventually zero. The group-cohomology side also vanishes, by Theorem \ref{thm_delta_stratified} on coherent subobjects and commutation of group cohomology with filtered direct limits. This proves the assertion for an isotypical ind-object, and the statement for finite direct sums follows.
\end{proof}

The restriction to finite direct sums in Proposition \ref{prop_elliptic}(2) cannot be removed in general.

\begin{example}[suggested by ChatGPT]\label{example_ordinary_elliptic_ind}\normalfont
Let $X/k$ be an ordinary elliptic curve. Choose
$0\neq Q_0\in\operatorname{Pic}^0(X)[p](k)$ and recursively
$Q_{m+1}\in\operatorname{Pic}^0(X)(k)$ with
$Q_{m+1}^{\otimes p}\simeq Q_m$; this is possible since multiplication by
$p$ on $\operatorname{Pic}^0(X)(k)$ is surjective. All $Q_m$ are
non-trivial. For $j\geq1$, define
\[
\mathcal L_j^{(n)}=
\begin{cases}
\mathcal O_X,&n<j,\\
Q_{n-j},&n\geq j.
\end{cases}
\]
Since $F^*Q\simeq Q^{\otimes p}$, together with
$Q_0^{\otimes p}\simeq\mathcal O_X$ and
$Q_{m+1}^{\otimes p}\simeq Q_m$, these terms carry the required Frobenius
transition isomorphisms and hence define a stratified line bundle
$\mathcal L_j$. Set
$\mathcal V_N=\bigoplus_{j=1}^N\mathcal L_j$ and
$\mathcal V=\bigoplus_{j\geq1}\mathcal L_j$.
For $n\geq N$, all summands of $\mathcal V_N^{(n)}$ are non-trivial
line bundles of degree zero, hence
\[
\mathrm H^q(X,\mathcal V_N^{(n)})=0,
\qquad q=0,1.
\]
Therefore Ogus' exact sequence \cite[Theorem~2.4]{Og75} gives
\[
\mathrm H^1_{\mathrm{str}}(X/k,\mathcal V_N)=0
\qquad\text{for all }N.
\]

For the infinite sum,
\[
A_n:=\mathrm H^0(X,\mathcal V^{(n)})=\bigoplus_{j>n}ke_j,
\]
where $e_j$ is the unit section of the $j$-th trivial summand. Regard
$(A_n,\varphi_n)$ as an inverse system of abelian groups; the transition
$\varphi_n:A_{n+1}\to A_n$ is Frobenius on the common coordinates. The
sequence $y_n=e_{n+1}$ defines a non-zero class in
${\varprojlim}^{1}_nA_n$: if $y_n=x_n-\varphi_n(x_{n+1})$, then for every
$m$
\[
x_0=e_1+\cdots+e_m+\varphi_{0,m}(x_m),
\]
which is impossible since $x_0\in A_0=\bigoplus_{j\geq1}ke_j$ has finite
support. Thus ${\varprojlim}^{1}_nA_n\neq0$, and
\cite[Theorem~2.4]{Og75} gives
$\mathrm H^1_{\mathrm{str}}(X/k,\mathcal V)\neq0$.
On the other hand, if $V_N=\mathcal V_N|_x$ and $V=\mathcal V|_x$, then
Theorem \ref{thm_delta_stratified} and the commutation of group-scheme
cohomology with filtered direct limits give
\[
\mathrm H^1(\pi(X),V)
=\varinjlim_N\mathrm H^1(\pi(X),V_N)=0.
\]
Hence $\delta^1$ is not bijective for $\mathcal V$.
\end{example}

\begin{rem}\label{rem_many_line_bundles}
The point of the preceding example is that, although each non-trivial
rank-one component eventually has vanishing cohomology, there is no uniform
Frobenius level from which this happens simultaneously for all components.
For every finite direct sum such a uniform level exists, which explains why
Proposition~\ref{prop_elliptic}(2) still applies in that case. For the infinite
direct sum above, the unbounded Frobenius delay produces a non-zero
${\varprojlim}^{1}$ term, and hence a stratified cohomology class which is not
detected at any finite stage.

This phenomenon is different from the counterexample constructed below for
curves of genus $2$. There the ind-stratified bundle lies entirely in the
unipotent isotypical component, and the failure of comparison is caused instead
by an unbounded order of Cartier nilpotency.
\end{rem}

\subsection{The order of Cartier nilpotency of a Frobenius invariant vector bundle} 
\label{sect_Cartier_nilpotency}
Next we show that the comparison map $\delta^i$ for
curves of genus at least 2 may be non-bijective. For a counter-example, we
construct an infinitely iterated extension of the trivial vector bundle which is
Frobenius invariant and investigate its stratified cohomology. For this, we
introduce the concept of the order of Cartier nilpotency of a Frobenius invariant
vector bundle, which seems interesting on its own. 

Our approach is based on the Cartier map. 
Let $\Omega_X$ denote the sheaf of  differential forms on $X$. 
Let $\mathcal B_X$ denote the image of $F_{X*}\mathcal O_X$ in $F_{X*}\Omega_X$ under the differential $F_{X*}d$. This is called the sheaf of locally exact forms. Then we have
exact sequences
$$0\longrightarrow \mathcal O_X\longrightarrow F_{X*}\mathcal O_X \longrightarrow
\mathcal B_X\longrightarrow 0$$
and
$$0\longrightarrow \mathcal B_X\longrightarrow 
F_{X*}\Omega_X\stackrel C\longrightarrow \Omega_X\longrightarrow 0,$$
where $C$ denotes the Cartier map \cite{Car57}. Roughly speaking, the differential 
forms that are not locally exact have the form $u^{p-1}du$ where $u$ is not a $p$-power
(in other words, $K$ is a separable extension of $k(u)$). Such a form is mapped to 
$du$ by the Cartier map $C$, other forms are mapped to zero by $C$.
Serre's duality yields a perfect pairing given by the residue map (cf. \cite{Ser58}):
$$\mathrm{H}^1(X,\mathcal O_X) \times \Omega(X) \longrightarrow k, \quad
(e,\mu)\longmapsto \langle e,\mu \rangle= \textrm{Res}(e_{ij} \mu):= \sum_{P\in |X|}\textrm{res}_P(f_i\mu).$$ 
Then $C$ is conjugate to $F^*_X$ with respect
to this pairing in the sense that
$\langle F^*_Xe,\mu\rangle =\langle e, C\mu\rangle ^p.$
In fact, we have, for any closed point $P$,
$\textrm{res}_P(\mu)=\textrm{res}_P(C\mu)^p.$

Consider the restriction of $C$ on the space $\Omega(X)$ of global 1-forms.
$C$ is $p^{-1}$-linear. Hence $\Omega(X)$ decomposes into the direct sum of a 
``$C$-bijective'' subspace and a ``$C$-nilpotent'' subspace.
This decomposition is compatible with the decomposition of 
$\mathrm H^1(X,\mathcal O_X)$ given above.

Let $\mathcal E$ be a vector bundle. 
Tensoring the second sequence with $\mathcal E$ and taking the cohomology
we get
$$\mathrm H^0(X,\mathcal E\otimes F_{X*} \Omega_X)\longrightarrow 
\mathrm H^0(X,\mathcal E\otimes  \Omega_X).$$
By the projection formula 
$$\mathcal E\otimes F_{X*}\Omega_X\cong F_{X*}(F^*\mathcal E\otimes \Omega_X) $$
and, as $F_X$ is affine,
$$\mathrm H^0(X,F_{X*}\mathcal H)\cong \mathrm H^0(X,\mathcal H) $$
for any quasi-coherent sheaf on $X$ (\cite[Ex.~8.1]{Har77}), hence we obtain a map
$$C_\mathcal{E}: \mathrm H^0(X,F^*\mathcal E\otimes \Omega_X)
\longrightarrow
\mathrm H^0(X,\mathcal E\otimes \Omega_X),$$
called \textit{the Cartier map} for $\mathcal E$. 
This map is conjugate to the Frobenius map on $\mathrm{H}^1(X,\mathcal E^\vee)\longrightarrow \mathrm{H}^1(X, F^*\mathcal{E}^\vee)$
by Serre's duality.

Let now $\mathcal E$ be Frobenius invariant. Fix an isomorphism
$\sigma: \mathcal E\longrightarrow F^*\mathcal E$.
We obtain a $p^{-1}$-linear map 
\begin{equation}\label{eq_Cartier_on_forms}
C_\sigma: \mathrm{H}^0(X,\Omega_X \otimes \mathcal E) \xrightarrow{\sigma} 
\mathrm{H}^0(X,\Omega_X \otimes F^* \mathcal E) \xrightarrow{C} 
\mathrm{H}^0(X,\Omega_X \otimes \mathcal E).\end{equation}

\begin{defn} \label{nilpotency} The order of Cartier nilpotency of the
Frobenius-equivariant bundle $(\mathcal E,\sigma)$ is the smallest integer $n$
such that $C_\sigma^n=0$ on the nilpotent
part $\mathrm{H}^0(X,\Omega_X\otimes \mathcal E)$. 
\end{defn}

The following result will provide a more comfortable way to compute the order of Cartier nilpotency.
\begin{prop} \label{proposition 5.6} Let $C^{n,0}$ denote the composition of $n$ Cartier morphisms:
$$C_i:=C_{F^{*(i-1)} \mathcal E}: \mathrm{H}^0(X,\Omega_X \otimes F^{*i} \mathcal E) \rightarrow \mathrm{H}^0(X,\Omega_X \otimes F^{*(i-1)}\mathcal E),$$ 
$$C^{n,0}:=C_1\circ C_2\circ\ldots\circ C_n:\mathrm{H}^0(X,\Omega_X \otimes F^{*n} \mathcal E) \rightarrow \mathrm{H}^0(X,\Omega_X \otimes \mathcal E).$$
The order of Cartier nilpotency of $(\mathcal E,\sigma)$ is equal to 
the smallest integer $n$ satisfying that $C^{n,0}$ vanishes on the nilpotent part of
$\mathrm{H}^0(X,\Omega_X \otimes F^{*n} \mathcal E)$.
\end{prop}
\begin{proof}  The construction of the Cartier map $C_\mathcal{E}$ is
functorial. In particular, the functoriality with respect to the map
$F^{*(i-1)}\sigma: F^{*(i-1)}\mathcal E\longrightarrow F^{*i}\mathcal E$ amounts
to the following commutative diagram: 
\begin{equation} \label{commutative diagram} \begin{tikzcd}
\mathrm{H}^0(X,\Omega_X \otimes F^{*i}\mathcal E) \arrow{r}{C_i} \arrow[swap]{d}{F^{*i}\sigma} 
& \mathrm{H}^0(X,\Omega_X \otimes F^{*(i-1)}\mathcal E) \arrow{d}{F^{*(i-1)}\sigma} \\%
\mathrm{H}^0(X,\Omega_X \otimes F^{*(i+1)}\mathcal E) \arrow{r}{C_{i+1}}& 
\mathrm{H}^0(X,\Omega_X \otimes F^{*i}\mathcal E),
\end{tikzcd}
\end{equation}
for any $i.$ Hence, one can check easily by induction that
$$(C\circ \sigma)^n= C_1\circ C_2\circ\ldots\circ C_n\circ F^{*(n-1)}\sigma\circ
\ldots\circ \sigma.$$
The claim follows. 
\end{proof}

\subsection{The construction of the example}\label{sec:construction-example}
Throughout this subsection, $X$ is a smooth projective curve over the
algebraically closed field $k$ of characteristic $p>0$, with
\begin{equation*}\tag{A}\label{assumption_genus_two}
g(X)=2,\qquad
\dim_k\mathrm H^1(X,\mathcal O_X)_\mathrm{ss}=1.
\end{equation*}
Thus $X$ has $p$-rank one; equivalently, the stable rank of Frobenius on
$\mathrm H^1(X,\mathcal O_X)$ is one.
The Frobenius--Cartier adjunction under Serre duality identifies the duals
of the semisimple and nilpotent summands of $\mathrm H^1(X,\mathcal O_X)$
with the corresponding summands of $\mathrm H^0(X,\Omega_X)$.
Both have dimension one, so we fix
\[
\mathrm H^0(X,\Omega_X)=k\eta\oplus k\omega,
\qquad C(\eta)=\eta,\qquad C(\omega)=0,
\]
where $\eta\ne0$ is chosen in the semisimple part and
$0\ne\omega\in\mathrm H^0(X,\Omega_X)_\mathrm{nil}$.
The normalization of $\eta$ is possible since $k$ is algebraically closed.

Our aim is to construct an ind-stratified bundle for which the comparison
map $\delta^1$ is not bijective.  We first choose
\[
0\ne e_2\in\mathrm H^1(X,\mathcal O_X)_\mathrm{ss},
\qquad F^*e_2=e_2.
\]
Such a class is unique up to multiplication by an element of
$\mathbb F_p^\times$; see \cite[\S14, Corollary, p.~143]{Mu70}.
The Artin--Schreier sequence identifies the Frobenius-fixed classes in
$\mathrm H^1(X,\mathcal O_X)$ with
$\mathrm H^1_{\mathrm{et}}(X,\mathbb F_p)$.  Since the $p$-rank of $X$ is
one, the maximal pro-$p$ quotient of $\pi_1^{\mathrm{et}}(X)$ is
$\mathbb Z_p$ \cite[Theorem~1.9]{Cre84}.  We therefore choose a quotient
\[
q:\pi_1^{\mathrm{et}}(X)\twoheadrightarrow\mathbb Z_p
\]
whose reduction modulo $p$ corresponds to $e_2$.

For $r\geq1$, put $V_r=\mathbb F_p[t]/(t^r)$ and let a topological generator of
$\mathbb Z_p$ act by multiplication by $1+t$.  This action is continuous,
since it factors through a finite cyclic $p$-group.  Let
$(\mathcal E_r,\sigma_r)$ be the corresponding finite stratified bundle,
where
\[
\sigma_r:\mathcal E_r\xrightarrow{\sim}F^*\mathcal E_r
\]
is the compatible Frobenius isomorphism supplied by this chosen
$\mathbb F_p$-local system.  In particular,
$\mathcal E_1=\mathcal O_X$, and the extension defining
$\mathcal E_2$ has class $e_2$.

We fix once and for all these compatible Frobenius isomorphisms
$\sigma_r:\mathcal E_r\xrightarrow{\sim}F^*\mathcal E_r$; the Cartier
operators below are defined relative to this fixed system.

The exact sequence
$0\to V_r\xrightarrow{t}V_{r+1}\to\mathbb F_p\to0$ becomes the extension
defining $\mathcal E_{r+1}$.  Its pushout along $V_r\to V_1$ is the extension
defining $\mathcal E_2$, whose underlying vector-bundle class is $e_2$;
hence every extension in the tower is nonsplit as a vector bundle.  The
associated cohomology sequences consequently give
\[
\mathrm H^0(X,\mathcal E_r)=k,
\qquad \dim_k\mathrm H^1(X,\mathcal E_r)_\mathrm{ss}=1,
\]
and the maps
\[
\mathrm H^1(X,\mathcal E_r)_\mathrm{ss}
\longrightarrow \mathrm H^1(X,\mathcal E_{r+1})_\mathrm{ss}
\]
are zero.  Thus $\mathcal E:=\varinjlim_r\mathcal E_r$ is the desired
candidate; the remaining point is to impose the Cartier condition that makes
its stratified cohomology large.

The following exact sequences are now part of the construction.
\begin{prop}\label{exact sequence of En}
For any positive integers $m$ and $n$, there is a Frobenius-equivariant
exact sequence of vector bundles over $X$:
\[
0 \longrightarrow \mathcal E_n \longrightarrow \mathcal E_{m+n}
\longrightarrow \mathcal E_m \longrightarrow 0.
\]
\end{prop}
\begin{proof}
On the $\mathbb Z_p$-modules used above, multiplication by $t^m$ and reduction
modulo $t^m$ give an exact sequence
\[
0\longrightarrow \mathbb F_p[t]/(t^n)
\xrightarrow{t^m}\mathbb F_p[t]/(t^{m+n})
\longrightarrow \mathbb F_p[t]/(t^m)\longrightarrow0.
\]
Both maps commute with multiplication by $1+t$.  Applying the exact functor
from these finite representations to finite stratified bundles gives the
asserted Frobenius-equivariant sequence.
\end{proof}


For $m,n\geq1$, write the exact sequence of Proposition
\ref{exact sequence of En} as
\[
0\longrightarrow\mathcal E_n
\xrightarrow{i_{n,m+n}}\mathcal E_{m+n}
\xrightarrow{\mathrm{pr}_{m+n,m}}\mathcal E_m
\longrightarrow0,
\]
and use the same notation for the induced maps on differential forms.
The Frobenius equivariance proved in Proposition
\ref{exact sequence of En} says explicitly that
\[
F^*i_{n,m+n}\circ\sigma_n
=\sigma_{m+n}\circ i_{n,m+n},\qquad
F^*\mathrm{pr}_{m+n,m}\circ\sigma_{m+n}
=\sigma_m\circ\mathrm{pr}_{m+n,m}.
\]
Set
\[
M_r:=\mathrm H^0(X,\Omega_X\otimes\mathcal E_r)_\mathrm{nil}
\quad(r\geq1),\qquad M_0:=0.
\]
Denote the resulting Cartier operator on $M_r$ by $C_r$.
In particular, $M_1=k\omega$; all inverse images below are taken inside $M_n$.
We isolate the additional hypothesis used in the counterexample:
\begin{equation*}\tag{C}\label{condition_Cartier}
C_2(M_2)=0,\qquad C_3(M_3)\subset i_{1,3}(M_1),\qquad
\overline C_3:M_3/i_{2,3}(M_2)\xrightarrow{\sim}M_1.
\end{equation*}
Here $\overline C_3$ is the $p^{-1}$-linear map induced by $C_3$, with
$i_{1,3}(M_1)$ identified with $M_1$. Naturality gives
$C_3\circ i_{2,3}=i_{2,3}\circ C_2=0$, so $C_3$ indeed factors through
$M_3/i_{2,3}(M_2)$. Equivalently, in addition to $C_2(M_2)=0$, a
nilpotent lift $s_3$ of $\omega$ satisfies
$C_3(s_3)=a\,i_{1,3}(\omega)$ for some $a\ne0$; the same relation then
holds for every such lift.
Condition \eqref{condition_Cartier} implies Cartier nilpotency order two
on $M_3$.
It will be verified in characteristic three in Section \ref{sect_genus_2}.

\begin{lem}\label{lemma 5.13}\leavevmode
\begin{itemize}
\item[i)] For all $m,n\geq1$, there is an exact sequence
\[
0\longrightarrow M_n\xrightarrow{i_{n,n+m}}M_{n+m}
\xrightarrow{\mathrm{pr}_{n+m,m}}M_m\longrightarrow0.
\]
In particular, $\dim M_n=n$, and every $M_n$ contains a lift of $\omega$
under $\mathrm{pr}_{n,1}$ (with $\mathrm{pr}_{1,1}=\mathrm{id}$).
\item[ii)] Assume \eqref{condition_Cartier}. For $n>2$ and
$s\in M_n$ with $\mathrm{pr}_{n,1}(s)=\omega$, there exist
$s'\in M_{n-2}$ and $a\in k^\times$ such that
\[
\mathrm{pr}_{n-2,1}(s')=\omega,
\qquad C_n(s)=a\,i_{n-2,n}(s').
\]
\end{itemize}
\end{lem}
\begin{proof}
For i), tensor the sequence in Proposition \ref{exact sequence of En}
with $\Omega_X$ and take cohomology. The connecting morphism is
Cartier-equivariant. By Serre duality, the Cartier operator on
$\mathrm H^1(X,\Omega_X\otimes\mathcal E_n)$ is dual to Frobenius on
$\mathrm H^0(X,\mathcal E_n^\vee)\simeq k$, hence is bijective. Thus the
target of the connecting morphism has trivial nilpotent part. Since the
semisimple--nilpotent decomposition of a finite-dimensional vector space
with semilinear endomorphism is functorial for equivariant maps, passing to
nilpotent parts gives the asserted short exact sequence. Its dimension
statement follows from $M_1=k\omega$.

For ii), put $s_3=\mathrm{pr}_{n,3}(s)$. The projection
$\mathrm{pr}_{n,3}:\mathcal E_n\to\mathcal E_3$ is compatible with the
chosen Frobenius isomorphisms. Naturality of the Cartier map for this
bundle morphism therefore gives the commutative diagram
\[
\begin{tikzcd}
M_n \arrow[r,"C_n"] \arrow[d,"\mathrm{pr}_{n,3}"']
& M_n \arrow[d,"\mathrm{pr}_{n,3}"]\\
M_3 \arrow[r,"C_3"'] & M_3.
\end{tikzcd}
\]
Thus condition \eqref{condition_Cartier} gives, concretely,
\[
\mathrm{pr}_{n,3}\bigl(C_n(s)\bigr)
=C_3\bigl(\mathrm{pr}_{n,3}(s)\bigr)
=C_3(s_3)=a\,i_{1,3}(\omega),\qquad a\ne0.
\]
Since $\mathrm{pr}_{n,2}=\mathrm{pr}_{3,2}\mathrm{pr}_{n,3}$ and
$\mathrm{pr}_{3,2}i_{1,3}=0$, it follows that
$\mathrm{pr}_{n,2}(C_n(s))=0$. Applying i) with $m=2$ and with its
index $n$ replaced by $n-2$, we obtain
$C_n(s)=i_{n-2,n}(u)$ for some $u\in M_{n-2}$.
Finally, the restriction of $\mathrm{pr}_{n,3}$ to the subbundle
$\mathcal E_{n-2}\subset\mathcal E_n$ is the composite
\[
\mathcal E_{n-2}\xrightarrow{\mathrm{pr}_{n-2,1}}\mathcal E_1
\xrightarrow{i_{1,3}}\mathcal E_3.
\]
Consequently
\[
\mathrm{pr}_{n,3}i_{n-2,n}
=i_{1,3}\mathrm{pr}_{n-2,1}.
\]
Comparison with the preceding display gives
$\mathrm{pr}_{n-2,1}(u)=a\omega$.
Taking $s'=a^{-1}u$ proves ii).
\end{proof}

\begin{prop}\label{Prop_Cartier_nilpotency}
Under \eqref{condition_Cartier}, the order of Cartier nilpotency of
$\mathcal E_n$ is at least $\lfloor(n+1)/2\rfloor$ for every $n\geq1$.
\end{prop}
\begin{proof}
Choose $s_n\in M_n$ mapping to $\omega$. Repeatedly applying
Lemma \ref{lemma 5.13} ii), with the index decreasing by two, gives
\[
C_n^j(s_n)=a_j\,i_{n-2j,n}(s_{n-2j}),\qquad a_j\ne0,
\quad 0\leq j\leq\left\lfloor\frac{n-1}{2}\right\rfloor,
\]
where every $s_{n-2j}$ maps to $\omega$; set $i_{n,n}=\mathrm{id}$.
The scalars remain non-zero under the $p^{-1}$-linear operator $C_n$.
At the last step the lift lies in $M_1$ or $M_2$, and is non-zero.
Hence $C_n^{\lfloor(n-1)/2\rfloor}(s_n)\ne0$, which proves the bound.
\end{proof}

\begin{prop}\label{prop-non-example}
Assume \eqref{assumption_genus_two} and \eqref{condition_Cartier}. Then
for the Frobenius invariant vector bundles $\mathcal E_r$ with inductive limit
$\mathcal E$ constructed above, and with $E:=\mathcal E|_x$, the comparison map
\[
\delta^1_\mathcal{E}: \mathrm{H}^1(\pi(X),E)=
\varinjlim_r \mathrm{H}^1_\mathrm{str}(X,\mathcal E_r)\longrightarrow
\mathrm{H}^1_\mathrm{str}(X,\mathcal E)
\]
is not bijective. In fact, the target space is infinite dimensional while the source space 
is zero. 
\end{prop} 
\begin{proof}
As observed above, the semisimple parts
$\mathrm H^1(X,\mathcal E_r)_\mathrm{ss}$ are one-dimensional and are
killed by the transition maps. Since
$\mathrm H^1_\mathrm{str}(X,\mathcal E_r)=
\mathrm H^1(X,\mathcal E_r)_\mathrm{ss}$, we obtain
\[
\varinjlim_r \mathrm H^1_\mathrm{str}(X,\mathcal E_r)
=\varinjlim_r \mathrm H^1(X,\mathcal E_r)_\mathrm{ss}=0.
\]

We now analyze the inverse system occurring in stratified cohomology.  Using
the compatible Frobenius isomorphisms $\sigma_r$ fixed in the construction,
represent $\mathcal E_r$ by the $F$-divided system
\[
\mathcal E_r^{(n)}=\mathcal E_r\qquad (n\geq 0),
\]
whose transition maps are $\sigma_r^{-1}$. Passing to the filtered direct limit over $r$, we obtain
\[
\mathcal E^{(n)}\simeq\mathcal E\qquad\text{for every }n,
\]
compatibly with the Frobenius transition maps.

Write
\[
\Phi_r:=\mathrm H^1(X,\sigma_r^{-1})\circ F_X^*
:\mathrm H^1(X,\mathcal E_r)\longrightarrow
\mathrm H^1(X,\mathcal E_r)
\]
for the $p$-linear cohomological Frobenius determined by $\sigma_r$, and put
\[
N_r:=\mathrm H^1(X,\mathcal E_r)_\mathrm{nil}.
\]
Thus the subscript $\mathrm{nil}$ refers to the nilpotent part for $\Phi_r$.
The kernel of $\mathrm H^1(X,\mathcal E_r)\to\mathrm H^1(X,\mathcal E_{r+1})$
is the semisimple line spanned by the extension class. Hence $N_r\subset N_{r+1}$;
Riemann--Roch and $h^0(X,\mathcal E_r)=1$ give $\dim N_r=r$. Set $N_0=0$. Since every transition map kills the semisimple line whereas its restriction
to $N_r$ is injective, the filtered direct limit identifies with
\[
\mathrm H^1(X,\mathcal E)=\varinjlim_r\mathrm H^1(X,\mathcal E_r)
=\bigcup_rN_r=:N,
\]
and $N$ is infinite dimensional.  Compatibility of the $\sigma_r$ induces a
$p$-linear endomorphism $\Phi:N\to N$ whose restriction to $N_r$ is
$\Phi_r$.

We next identify $\Phi$ on the graded pieces of $N$. Put
$V_q=\mathbb F_p[t]/(t^q)$. On the dual $\mathbb Z_p$-module, the element $t$
acts as
\[
(1+t)^{-1}-1=-\frac{t}{1+t}=-t+t^2-\cdots.
\]
Since this polynomial has a non-zero linear term, its action is again one
nilpotent Jordan block of length $q$.
We may therefore choose a $\mathbb Z_p$-equivariant self-duality
$\beta_q:V_q\xrightarrow{\sim}V_q^\vee$. Being a morphism of
$\mathbb F_p$-local systems, it induces a Frobenius-compatible self-duality
$\beta_q:\mathcal E_q\xrightarrow{\sim}\mathcal E_q^\vee$.

We use this self-duality only in the fixed rank $q$; no compatibility between
the choices for different $q$ is needed. It exchanges the inclusions and
projections in the tower. More precisely, when $q=m+n$, the dual of
\[
0\longrightarrow\mathcal E_n\longrightarrow\mathcal E_q
\longrightarrow\mathcal E_m\longrightarrow0
\]
is identified via $\beta_q$, and the induced identifications of the end terms,
with the reversed sequence
\[
0\longrightarrow\mathcal E_m\longrightarrow\mathcal E_q
\longrightarrow\mathcal E_n\longrightarrow0.
\]
Indeed, each Jordan block has a unique submodule of every dimension.

Use $\beta_q$ in Serre duality. Since $\Phi_q$ and $C_q$ are adjoint, the
nilpotent parts are in perfect duality:
\[
\langle\ ,\ \rangle_q:N_q\times M_q\longrightarrow k.
\]
We regard $N_j$ and $M_j$ as subspaces of $N_q$ and $M_q$, respectively,
using the maps in the tower.  Naturality of Serre duality, together with
the reversal of the two flags under the chosen self-duality, gives
\[
(N_j)^\perp=M_{q-j}\qquad(0\leq j\leq q).
\]
Consequently
\[
(N_j/N_{j-1})^\vee\simeq M_{q-j+1}/M_{q-j}\qquad(1\leq j\leq q).
\]
Now fix $r\geq1$ and take $q=r+2$.  For $j=q$ and $j=r$, respectively,
this gives
\[
(N_{r+2}/N_{r+1})^\vee\simeq M_1,\qquad
(N_r/N_{r-1})^\vee\simeq M_3/M_2.
\]
Since cohomological Frobenius and Cartier are adjoint, the transpose of the
map induced by $\Phi$ is precisely
$\overline C_3:M_3/M_2\xrightarrow{\sim}M_1$, which is an isomorphism by
\eqref{condition_Cartier}. Therefore
\[
\Phi:N_{r+2}/N_{r+1}\xrightarrow{\sim}N_r/N_{r-1}
\qquad(r\geq1).
\]
These graded isomorphisms imply $\Phi:N\twoheadrightarrow N$. To see this,
argue by induction on $r$. Given $v\in N_r$, choose $w_0\in N_{r+2}$ such
that $\Phi(w_0)\equiv v\pmod{N_{r-1}}$ by the preceding isomorphism. The
difference $v-\Phi(w_0)$ lies in $N_{r-1}$
and therefore has a preimage by the induction hypothesis.

By the exact sequence recalled in Subsection \ref{sect_Ogus}, we have
\[
0\longrightarrow {\varprojlim}^{1}_n\mathrm H^0(X,\mathcal E^{(n)})
\longrightarrow \mathrm H^1_\mathrm{str}(X,\mathcal E)
\longrightarrow \varprojlim_n\mathrm H^1(X,\mathcal E^{(n)})
\longrightarrow 0.
\]
For every $n$, we have $\mathrm H^0(X,\mathcal E^{(n)})\simeq k$, and the transition maps are bijective $p$-linear maps on $k$. Therefore
${\varprojlim}^{1}_n\mathrm H^0(X,\mathcal E^{(n)})=0$, and we obtain a natural isomorphism
\[
\mathrm H^1_\mathrm{str}(X,\mathcal E)\xrightarrow{\sim}
\varprojlim_n \mathrm H^1(X,\mathcal E^{(n)}).
\]
Under the compatible identifications above, the inverse system on the right is
\[
N\xleftarrow{\Phi}N\xleftarrow{\Phi}N\xleftarrow{\Phi}\cdots .
\]
Since $\Phi:N\twoheadrightarrow N$, the projection
\[
\varprojlim_\Phi N\longrightarrow N,\qquad (v_0,v_1,\ldots)\longmapsto v_0,
\]
is surjective. Hence $\varprojlim_\Phi N$, and therefore
$\mathrm H^1_\mathrm{str}(X,\mathcal E)$, is infinite dimensional.
\end{proof}

\subsection{Explicit computation for hyperelliptic curves of genus 2 over a field of characteristic 3}
\label{sect_genus_2}
We now verify \eqref{condition_Cartier} on explicit curves in
characteristic three. The characteristic is used only in this realization;
the arguments of Section \ref{sec:construction-example} require
\eqref{assumption_genus_two} and \eqref{condition_Cartier}.
Let $X$ be a hyperelliptic curve of genus 2 with two marked Weierstrass points, and 
$\text{char}(k)=3.$ By fixing two points, such curve can be defined by the following affine covering: $X = U_0 \cup U_1,$ where
\begin{align*}
	U_0 = X \setminus \{\infty\}: y^2 = a_5x^5+a_4x^4+a_3x^3+a_2x^2+a_1x; \\
	U_1 = X \setminus \{0\}: w^2 = a_1v^5+a_2v^4+a_3v^3+a_4v^2+a_5v,
\end{align*}
and on the intersection $U_{01}:= U_0 \cap U_1,$ we have the following relations: $y =  w v^{-3}$
and $x = v^{-1}.$ 
We first provide the necessary and sufficient conditions for $X$ to have Hasse-Witt rank $1$. 

\begin{lem} \label{Hasse-Witt rank one}
Let $X$ be the curve obtained by patching the two affine curves
\begin{align*}
	U_0 = X \setminus \{\infty\}: y^2 = a_5x^5+a_4x^4+a_3x^3+a_2x^2+a_1x; \\
	U_1 = X \setminus \{0\}: w^2 = a_1v^5+a_2v^4+a_3v^3+a_4v^2+a_5v,
\end{align*}
subject to the relations: $y =  w v^{-3}$ and $x = v^{-1}$. 
Then  $X$ is smooth (hence has genus 2 and is hyperelliptic)
and has Hasse-Witt rank one  if and only if 
$$\begin{cases} a_4a_2 = a_5a_1 \\ a_1a_5(a_2^{4} + a_1^{3}a_5) \neq 0\\
a_3(a_2a_5^2+a_4^3-a_3a_4a_5)\neq 0. \end{cases}$$
\end{lem}  
\noindent The proof will be given in the Appendix \ref{proof_1}. \hfill $\Box$

\begin{lem}[The vector bundle $\mathcal E_2$] \label{lem_E2}
Assume the conditions of Lemma \ref{Hasse-Witt rank one}.
The cohomology class of the 1-cocycle $e_{01} \in \mathcal{O}(U_{01})$  
\begin{equation} \label{e} e_{01} = ( a + bx)x^{-2}y = a.\frac{w}{v} +b.\frac{y}{x}, \end{equation}
is invariant under Frobenius if $(a,b)\ne(0,0)$ satisfies:
\begin{equation} \label{e01} 
\begin{bmatrix}
	a_4 & a_1 \\ a_5 & a_2
\end{bmatrix} \cdot \begin{bmatrix}
a^3 \\ b^3
\end{bmatrix} = \begin{bmatrix}
a \\ b
\end{bmatrix}.
\end{equation}\end{lem}
\noindent The proof will be given in the Appendix \ref{proof_2}. \hfill $\Box$

\begin{lem}\label{lemma 3.7}
For the bundle $\mathcal E_2$ of Lemma \ref{lem_E2}, the induced
Cartier operator vanishes on $M_2$; its nilpotency order is one.  
\end{lem}
\noindent The proof will be given in the Appendix \ref{proof_3}. \hfill $\Box$

\begin{example}
	If $a_i = -1$ for all $1 \leq i \leq 5,$ then we can choose
	\begin{align*}
		\omega  = \frac{(x-1)dx}{y}; \quad  \eta  = \frac{(x+1)dx}{y}\\
		e_{01} = (1+x)\frac{y}{x^2} = \frac{w}{v} + \frac{w}{v^2}.
	\end{align*}
\end{example}

\begin{lem}[The vector bundle $\mathcal E_3$]\label{lem_E3}
The bundle $\mathcal E_3$ is given by a $GL_3(\mathcal{O}(U_{01}))$-cocycle of the form:
\begin{equation*}
A_3= \begin{bmatrix}
1 & e_{01} & f_{01} \\  0 & 1 & e_{01} \\ 0 & 0 & 1
\end{bmatrix}.
\end{equation*}with $f_{01} = -e_{01}^2$. \end{lem}
\noindent The proof will be given in the Appendix \ref{proof_4}. \hfill $\Box$

After the change of trivialization in Appendix \ref{proof_4}, we use
the cocycle \eqref{cocycle3}, transporting the Frobenius isomorphism with it.
Then each element in 
$\mathrm{H}^0(X,\Omega_X \otimes\mathcal  E_3)$ is a pair of vectors   
of local sections subject to the relation
\begin{align*}
& \begin{bmatrix} 
{}^3s_0 \\ {}^2s_0 \\ {}^1s_0 \end{bmatrix} =\begin{bmatrix}
1 & e_{01} & 0 \\  0 & 1 & e_{01} \\ 0 & 0 & 1
\end{bmatrix}. \begin{bmatrix} {}^3s_1 \\ {}^2s_1 \\ {}^1s_1 \end{bmatrix}
\end{align*}
where ${}^j\!s_i\in\Omega(U_i)$ for $i=0,1$ and $j=1,2,3$, that is 
$$ \begin{cases} {}^3s_0 - {}^3s_1 = e_{01}. {}^2s_1 \\ {}^2s_0 - {}^2s_1 = e_{01}. {}^1s_1 \\ {}^1s_0 - {}^1s_1 = 0. \end{cases} $$
Since ${}^1s_0 - {}^1s_1 = 0,$ these local sections (of the sheaf $\Omega_X$) glue together to obtain a global section $s$ in $\Omega(X).$

Our finding here is the following.
\begin{lem}\label{lem_C_E3}
Assume the conditions of Lemma \ref{Hasse-Witt rank one} and use the
bundles $\mathcal E_2,\mathcal E_3$ constructed above.
Then condition \eqref{condition_Cartier} holds. In particular, the
Cartier nilpotency order of $\mathcal E_3$ is two.
\end{lem}
\noindent The proof is given in Appendix \ref{proof_5}. \hfill $\Box$

Together with Propositions \ref{Prop_Cartier_nilpotency} and
\ref{prop-non-example}, this verifies the counterexample in characteristic three.

\begin{example}\label{example 5.11}
\begin{itemize}
\item[(i)] The curve
$y^2=-x^5-x^4-x^3-x^2-x$, defined over $\mathbb F_3$ and base-changed
to $k$, satisfies Lemma \ref{Hasse-Witt rank one}, hence
\eqref{condition_Cartier}.
\item[(ii)] The smooth curve
$y^2=-x^5+x^4+x^3-x^2+x$ has $p$-rank zero:
here $a_2^4+a_1^3a_5=0$, and its non-zero Cartier operator satisfies $C^2=0$.
Thus it does not satisfy \eqref{assumption_genus_two} and cannot be used
in this construction.
\end{itemize}
\end{example}

\begin{appendix}
\section{}  
\subsection{Proof of Lemma \ref{Hasse-Witt rank one}}\label{proof_1}
$X$ is smooth if and only if the degree-five polynomial $f(x)$ has no repeated roots. This is the same as saying that $a_1a_5 \neq 0$ and the discriminant of the quartic polynomial $a_5x^4 + a_4x^3+a_3x^2+a_2x+a_1,$ 
\begin{align*}
\Delta := {}& a_5^3a_1^3+a_5^2a_3^2a_1^2+a_5a_4a_3^2a_2a_1+a_5a_3^4a_1\\
&{}-a_5a_3^3a_2^2-a_4^3a_2^3-a_4^2a_3^3a_1+a_4^2a_3^2a_2^2\ne0. 
\end{align*} 

The space of global differential 1-forms $\Omega(X)$
is spanned by the holomorphic differential 
$$ \omega_1 = y^{-1}dx \quad \text{ and } \quad  \omega_2= y^{-1}xdx.$$
The Cartier operator in this basis is given by:
$$C = \begin{bmatrix}
	a_2^{1/3} & a_1^{1/3} \\ a_5^{1/3} & a_4^{1/3}
\end{bmatrix}.
$$
For $X$ to have Hasse-Witt rank one, it is necessary that $\det(C)=0$, that is
$$a_1a_5 = a_2a_4.$$
Under this equality, the condition $\Delta\neq 0$ can be expressed as
$$a_3(a_2a_5^2+a_4^3-a_3a_4a_5)\neq 0.$$

Since $C$ has rank one, its stable rank is one exactly when $C^2\ne0$.
In that case there is a basis $\{\omega,\eta\}$ of $\Omega(X)$
with $C(\omega)=0$ and $C(\eta)=\eta$. 
We can choose
$$\omega = a_4 \omega_1 - a_5 \omega_2.$$
If we present $\eta$ in the form $\eta = \lambda_1^3 \omega_1 + \lambda_2^3 \omega_2$
then
$$C(\eta) = C(\lambda_1^3 \omega_1 + \lambda_2^3 \omega_2) = \lambda_1C(\omega_1) + \lambda_2C(\omega_2).$$
So $\lambda_1, \lambda_2$ need to satisfy:
\begin{align*} 
& \begin{bmatrix}
	a_2^{1/3} & a_1^{1/3} \\ a_5^{1/3} & a_4^{1/3}
\end{bmatrix} \cdot \begin{bmatrix}
\lambda_1 \\ \lambda_2
\end{bmatrix} = \begin{bmatrix}
\lambda_1^3 \\ \lambda_2^3
\end{bmatrix} \end{align*}
Elementary computation shows that the existence of $\eta \neq 0$ is equivalent to   
$$a_2^{4} + a_1^{3}a_5 \neq 0.$$
 
\subsection{Proof of Lemma \ref{lem_E2}}\label{proof_2}
One needs to check 
$e_{01}^3-e_{01} = 0$
 in $\mathrm{H}^1(X,\mathcal{O}_X).$ 
Indeed, we have
\begin{align*}
	e_{01}^3 &  = a^3\frac{w^3}{v^3} + b^3\frac{y^3}{x^3} \\
	& =  a^3\frac{w(a_1v^5+a_2v^4+a_3v^3+a_4v^2+a_5v^1)}{v^3} + b^3\frac{y(a_5x^5+a_4x^4+a_3x^3+a_2x^2+a_1x^1)}{x^3} \\
	& = a^3w(a_1v^2+a_2v+a_3) +b^3y(a_5x^2+a_4x+a_3) + \big(a^3a_4\frac{w}{v} + a^3a_5\frac{w}{v^2} + b^3a_2\frac{y}{x} + b^3a_1\frac{y}{x^2} \big) \\ 
	& = f_1+f_0 + \big(a^3a_4\frac{w}{v} + a^3a_5\frac{y}{x} + b^3a_2\frac{y}{x} + b^3a_1\frac{w}{v} \big) \\
	& = f_1+f_0 + \big(\frac{w}{v}(a^3a_4+b^3a_1) +\frac{y}{x}(a^3a_5+b^3a_2) \big) \\
	& = f_1+f_0 + \big(\frac{w}{v}.a +\frac{y}{x}.b \big) \quad \text{by (\ref{e01})} \\
	& = f_1 + f_0 + e_{01}
	\end{align*} where $f_0 = b^3y(a_5x^2+a_4x+a_3) \in \mathcal{O}(U_0),$ and $f_1 = a^3w(a_1v^2+a_2v+a_3) \in \mathcal{O}(U_1).$ 
 
\subsection{Proof of Lemma \ref{lemma 3.7}}\label{proof_3}
We now compute the image of $\begin{bmatrix} {}^2s_i \\ \omega \end{bmatrix} \in \mathrm{H}^0(X,\Omega_X \otimes F^{*}\mathcal E_2) $ under $C,$ where recall that $\omega = a_4 \omega_1 - a_5 \omega_2.$
According to \eqref{e}, $e_{01}^3 = \frac{y^3}{x^6}(a^3+b^3x^3)$. Hence 
\begin{align*}
e_{01}^3.\omega & = \frac{y^3}{x^6}(a^3+b^3x^3) . \omega = \frac{y^2}{x^6}(a^3+b^3x^3) . (a_4-a_5x)dx \\
& = \frac{(b^3x^3+a^3)(a_4-a_5x)(a_5x^4+a_4x^3+a_3x^2+a_2x+a_1)}{x^5}dx \\
& = \big(-b^3a_5^2x^3+b^3(a_4^2-a_5a_3)x+b^3(a_4a_3-a_5a_2)-a^3a_5^2\big)dx\\
&\qquad - \big(b^3a_4a_1 +a^3(a_4^2-a_5a_3) + a^3(a_4a_3-a_5a_2)v +a^3a_4a_1v^3\big)dv.
\end{align*} 
Consequently 
\begin{align*}
{}^2s_0 - {}^2s_1 & = e_{01}^3.\omega \\
& = \big(-b^3a_5^2x^3+b^3(a_4^2-a_5a_3)x+b^3(a_4a_3-a_5a_2)-a^3a_5^2\big)dx\\
&\qquad  -\big(b^3a_4a_1 +a^3(a_4^2-a_5a_3) + a^3(a_4a_3-a_5a_2)v +a^3a_4a_1v^3\big)dv.
\end{align*}
So, up to an element of $\Omega(X),$ we have that
\begin{align}\label{eq_2-cocycle0}
{}^2s_0 & = \big(-b^3a_5^2x^3+b^3(a_4^2-a_5a_3)x+b^3(a_4a_3-a_5a_2)-a^3a_5^2\big)dx \\
\label{eq_2-cocycle1}
{}^2s_1 & = \big(b^3a_4a_1 +a^3(a_4^2-a_5a_3) +a^3(a_4a_3-a_5a_2)v +a^3a_4a_1v^3\big)dv.
\end{align}
With these choices, $C({}^2s_0)=C({}^2s_1)=0$.
This lift of $\omega$, together with the section $[\omega,0]^t$,
gives two independent elements in the kernel of the raw Cartier map.
The section $[\eta,0]^t$ is not in its kernel. Since
$h^0(X,\Omega_X\otimes\mathcal E_2)=3$ and $\dim M_2=2$,
the kernel is exactly the nilpotent part (transported by the chosen
Frobenius isomorphism). Hence $C_2(M_2)=0$.

\subsection{Proof of Lemma \ref{lem_E3}}\label{proof_4}
We have 
$$e_{01}^3-e_{01} = f_0+f_1$$ with 
$f_0 = b^3y(a_5x^2+a_4x+a_3) \in \mathcal{O}(U_0)$ and 
$f_1 = a^3w(a_1v^2+a_2v+a_3) \in \mathcal{O}(U_1).$ 
Hence
\begin{align*}
  \begin{bmatrix}
1 & f_0 & -f_0^2 \\  0 & 1 & f_0 \\ 0 & 0 & 1
\end{bmatrix} . \begin{bmatrix}
1 & e_{01} & -e_{01}^2 \\  0 & 1 & e_{01} \\ 0 & 0 & 1
\end{bmatrix}\cdot \begin{bmatrix}
1 & f_1 & -f_{1}^2 \\  0 & 1 & f_{1} \\ 0 & 0 & 1
\end{bmatrix}  
%
&
= \begin{bmatrix}
1 & e_{01}^3 & -(f_1+f_0)^2 - 2e_{01}(f_0+f_1) - e_{01}^2 \\  0 & 1 & e_{01}^3 \\ 0 & 0 & 1
\end{bmatrix} \\
&=  \begin{bmatrix}
1 & e_{01}^3 & -e_{01}^6 \\  0 & 1 & e_{01}^3 \\ 0 & 0 & 1
\end{bmatrix}. 
\end{align*}
To pass to a simpler cocycle, write
\[
e_{01}^2=\frac{(a+bx)^2f(x)}{x^4}=g_0-g_1,
\qquad g_i\in\mathcal O(U_i),
\]
by splitting the Laurent polynomial into its non-negative and negative
powers of $x$. If $E_{13}$ denotes the elementary matrix, then
\[
(I+g_0E_{13})A_3(I+g_1E_{13})^{-1}
=\begin{bmatrix}1&e_{01}&0\\0&1&e_{01}\\0&0&1\end{bmatrix}.
\]
Thus we may use the equivalent cocycle
\begin{equation}\label{cocycle3}
\begin{bmatrix}1&e_{01}&0\\0&1&e_{01}\\0&0&1\end{bmatrix},
\end{equation}
with its transported Frobenius isomorphism. The changes of trivialization
preserve the flag and the quotient coordinates.

\subsection{Proof of Lemma \ref{lem_C_E3}}\label{proof_5}
We use the cocycle \eqref{cocycle3}. Thus a section of
$\Omega_X\otimes F^*\mathcal E_3$ with last coordinate $\omega$ satisfies
\[
{}^2s_0-{}^2s_1=e_{01}^3\omega,
\qquad {}^3s_0-{}^3s_1=e_{01}^3{}^2s_1.
\]
Choose ${}^2s_i$ as in \eqref{eq_2-cocycle0}--\eqref{eq_2-cocycle1},
and put
\[
\begin{aligned}
c_0&=b^3a_4a_1+a^3(a_4^2-a_5a_3),\\
c_1&=a^3(a_4a_3-a_5a_2),\qquad c_3=a^3a_4a_1.
\end{aligned}
\]
Then
\[
e_{01}^3{}^2s_1=
\frac{(a_1v^4+a_2v^3+a_3v^2+a_4v+a_5)^2}{wv^4}
(a^3v^3+b^3)(c_0+c_1v+c_3v^3)\,dv.
\]
Separating the negative powers of $v$ and using
$v^{-j}dv/w=-x^{j+1}dx/y$ gives
\[
e_{01}^3{}^2s_1=\tau_1-P(x)\frac{dx}{y},
\qquad \tau_1\in\Omega(U_1),\qquad
P(x)=A_2x^2+A_3x^3+A_4x^4+A_5x^5,
\]
where
\[
\begin{aligned}
A_2&=a_5^2(a^3c_0+b^3c_3)
       +b^3c_1(a_4^2-a_3a_5)-b^3c_0(a_2a_5+a_3a_4),\\
A_3&=b^3\bigl(-c_1a_4a_5+c_0(a_4^2-a_3a_5)\bigr),\\
A_4&=b^3(c_1a_5^2-c_0a_4a_5),\qquad A_5=b^3c_0a_5^2.
\end{aligned}
\]
Consequently the correct local choices are
\[
{}^3s_0=-P(x)\frac{dx}{y},\qquad {}^3s_1=-\tau_1.
\]
They may both be changed by the same global form $\gamma\in\Omega(X)$.
Since $C({}^2s_i)=C(\omega)=0$, their Cartier images glue to a global form.
The only powers in $f(x)P(x)$ that contribute to Cartier are $x^5$ and $x^8$;
their coefficients are
\[
R_5=a_3A_2+a_2A_3+a_1A_4,\qquad
R_8=a_5A_3+a_4A_4+a_3A_5.
\]
Equation \eqref{e01} and $a_1a_5=a_2a_4$ give
$(a,b)=t(a_4,a_5)$ for some $t\ne0$.
Substitution in these expressions yields
\[
R_8=0,\qquad
R_5=t^6a_3a_4^2a_5^2K,
\quad K=(a_5^2a_2+a_4^3)(a_5^2a_2+a_4^3-a_3a_4a_5).
\]
Both factors of $K$ are non-zero by Lemma \ref{Hasse-Witt rank one}:
for the first, use
$a_2^3(a_5^2a_2+a_4^3)=a_5^2(a_2^4+a_1^3a_5)$.
Thus $R_5\ne0$, and the raw Cartier map gives
\begin{equation}\label{C3s}
C({}^3s_0)=C({}^3s_1)
=-R_5^{1/3}\omega_2+C(\gamma),
\qquad \omega_2=\frac{x\,dx}{y}.
\end{equation}
Here $C(\gamma)\in k\eta$, whereas $\omega_2\notin k\eta$:
the image line of the Cartier matrix in Appendix \ref{proof_1}
has both coordinates non-zero. Hence \eqref{C3s} has a non-zero
nilpotent component.

Finally take the nilpotent component of the constructed section in
$\mathrm H^0(X,\Omega_X\otimes F^*\mathcal E_3)$, using the transported
Frobenius isomorphism. Its projection to $F^*\mathcal E_2$ is unchanged,
since that projection lies in the Cartier kernel by Appendix \ref{proof_3}.
Its image under the raw Cartier map is therefore a non-zero multiple of
$[\omega,0,0]^t$. Transporting back gives a lift $s_3\in M_3$ of $\omega$
with $C_3(s_3)=a\,i_{1,3}(\omega)$, $a\ne0$.
Together with $C_2(M_2)=0$ and Lemma \ref{lemma 5.13} i), this proves
\eqref{condition_Cartier}.

\subsection{An example with vanishing first stratified cohomology}
\label{proof_non-example} 
We give an example, on a hyperelliptic curve in characteristic \(3\), of an
\(F\)-periodic line bundle whose first stratified cohomology group vanishes.

Let $X$ be any smooth curve given by the equations of Section \ref{sect_genus_2};
we do not impose Hasse--Witt rank one in this subsection.
Let $\mathcal{L}:=\mathcal O(\infty-O)$ 
be the line bundle corresponding to the Weierstrass divisor $\infty-O,$ 
where $O$ is $(x,y)=(0,0)$ on $U_0$ and $\infty$ is $(v,w)=(0,0)$ on $U_1$.
Since $2\cdot\infty - 2\cdot O$ is the divisor of $x^{-1}$, $\mathcal L^{\otimes 2}$ is
trivial.
Hence $\mathcal{L}^3 \cong \mathcal{L}$, that is $\mathcal{L}$ is an $F$-periodic 
line bundle. Thus, $\mathcal L$ determines a stratified line bundle
which will also be denoted by $\mathcal{L}$. 

To compute $\mathrm H^1_\mathrm{str}(X,\mathcal L)$, we study Frobenius on
$\mathrm H^1(X,\mathcal L)$, or equivalently Cartier on
\[
\mathrm{H}^0(X,\Omega_X\otimes\mathcal L^\vee)=
\mathrm{H}^0(X,\Omega_X(O - \infty)).
\]

$\mathcal L$ is non-trivial: otherwise a rational function with divisor
$\infty-O$ would give a degree-one map $X\to\mathbb P^1$, contrary to $g(X)=2$.
Riemann--Roch therefore gives $h^0(X,\Omega_X(O-\infty))=1$, so there is a 
unique (up to scalars)  differential form which has a pole at $O$ of order at most
$1$, vanishes at $\infty$ and is regular elsewhere. 
Using the notations of Section \ref{sect_genus_2}, the form is:
\begin{equation*}
\omega_1 = \frac{dx}{y},
\end{equation*} 
and the Cartier map on $\mathrm{H}^0(X,\Omega(O - \infty))$ is the following composition
(cf. \eqref{eq_Cartier_on_forms}):
\begin{equation}
C_\sigma: \mathrm{H}^0(X,\Omega_X(O - \infty)) \xrightarrow{\cdot x^{-1}} 
\mathrm{H}^0(X,\Omega_X(3\cdot(O - \infty))) \xrightarrow{C} \mathrm{H}^0(X,\Omega_X(O - \infty)).
\end{equation}
The image of $\omega_1$ under this  map is 
\begin{align*}
C_\sigma(\omega_1) = C\left(\frac{y^2x^2dx}{x^3y^3}\right) & = \frac{1}{xy}
\cdot C((a_5x^5+a_4x^4+a_3x^3+a_2x^2+a_1x)x^2dx) \\
& = \frac{1}{y}\cdot a_3^{1/3}dx = a_3^{1/3}\cdot\omega_1.
\end{align*}
We conclude that $\mathrm H^1_\mathrm{str}(X,\mathcal L)=0$ if and only if the 
coefficient $a_3$ is 0. Notice that, according to the proof of \ref{Hasse-Witt rank one}, if
$a_3=0$ then $X$ has to be ordinary. Indeed, if $a_3=0$, 
$\Delta=(a_5a_1-a_4a_2)^3$. Thus, if $X$ is smooth, its Cartier matrix
is invertible and $X$ is ordinary. For example, $y^2=x^5+x$ over
$\overline{\mathbb F}_3$ is smooth and has $a_3=0$, giving the claimed example. 

\end{appendix} 

\section*{Acknowledgments}
We thank H\'el\`ene Esnault, Michel Brion, Joao Pedro dos Santos and Pham Thanh Tam for engaging in insightful discussions. 
The research of Ph\`ung H\^o Hai and D\`ao Van Thinh was supported by the Vietnam Academy of Science and Technology under grant numbers  CBCLCA.01/25-27 and CTTH00.02/25-26. The work of V\~o Qu\^oc Bao was supported  by the Vingroup Innovation Foundation under grant number VINIF.2023.TS.011 and by the ICTP/IAEA Sandwich Training Educational Programme (STEP).

\end{document}